\documentclass[12pt]{article}
\usepackage{amsmath,amssymb,amsbsy,url}
\usepackage{graphicx}
\usepackage{multirow}
\usepackage{algorithm}
\usepackage{xcolor}
\usepackage{tikz}
\usepackage{tkz-graph}
\usepackage[letterpaper, left=2.5cm, right=2.5cm, top=2.5cm,bottom=2.5cm,dvips]{geometry}

\newtheorem{proposition}{Proposition}[section]
\newtheorem{theorem}[proposition]{Theorem}
\newtheorem{lemma}[proposition]{Lemma}

\newtheorem{corollary}[proposition]{Corollary}

\newtheorem{definition}[proposition]{Definition}
\newtheorem{question}[proposition]{Question}

\newtheorem{problem}[proposition]{Problem}

\newcommand{\qed}{\hfill \rule{.1in}{.1in}}

\makeatletter
\def\imod#1{\allowbreak\mkern10mu({\operator@font mod}\,\,#1)}
\makeatother

\definecolor{dgreen}{rgb}{0.0, 0.5, 0.0}

\begin{document}

\title{On star b-chromatic number of a graph\thanks{The first author was partially suppored by Slovenian Research and Innovation Agency by program Telematics No. P2-0065. The third author was partially supported by the Slovenian Research and Innovation Agency by program No. P1-0297.}}
\author{Dragana Bo\v{z}ovi\'c$^{(1)}$, Da\v sa Mesari\v c \v Stesl$^{(2)}$, Iztok Peterin$^{(1,3)}$  \\
$^{(1)}$ {\small Faculty of Electrical Engineering and Computer Science}\\
{\small University of Maribor,} {\small Koro\v{s}ka cesta 46, 2000 Maribor,
Slovenia.} \\
$^{(2)}${\small Faculty of Computer and Information Science} \\
{\small University of Ljubljana,} {\small Ve\v{c}na pot 113, 1000 Ljubljana,
Slovenia.} \\
$^{(3)}${\small Institute of Mathematics, Physics and Mechanics}\\
{\small Jadranska ulica 19, 1000 Ljubljana, Slovenia.} \\
{\small \texttt{e-mails:}\textit{dragana.bozovic\@@um.si, iztok.peterin\@@um.si, dasa.stesl\@@fri.uni-lj.si}}
}
\date{}
\maketitle

\begin{abstract}
A star coloring of a graph $G$ is a proper coloring where vertices of every two color classes induce a forest of stars. A strict partial order is defined on the set of all star colorings of $G$. We introduce the star b-chromatic number $S_b(G)$, analogous to the b-chromatic number, as the maximum number of colors in a minimum element of the mentioned order. We present several combinatorial properties of $S_b(G)$, compute the exact value for $S_b(G)$ for several known families and compare $S_b(G)$ with several invariants naturally connected to $S_b(G)$.   
\end{abstract}

\noindent \textbf{Keywords}: star b-chromatic number; star b-coloring; b-coloring;  \medskip

\noindent \textbf{AMS subject classification (2020)}: 05C15

\section{Introduction}\label{sec_Intro}

The b-chromatic number $\varphi (G)$ of a graph $G$ was introduced in 1999 by Irving and Manlove \cite{IrMa} as a kind of a dual to chromatic number $\chi (G)$. They introduced a strict partial order on all proper colorings of $G$ where  $\chi(G)$ is the minimum number of colors and $\varphi(G)$ the maximum number of colors of a minimal element of the mentioned order. The topic gain quite some attention during the years, see the survey \cite{JaPe} and the references in there. In almost every publication on b-chromatic number an alternative definition of $\varphi(G)$ is used described already in \cite{IrMa} and is the maximum number of colors in a proper coloring where each color class contains a vertex with all colors in its closed neighborhood. One of the exceptions is the recent paper of Anholcer et al. \cite{AnCP} where the original definition of the b-chromatic number was used on the set of all acyclic colorings of a graph $G$ to introduce the acyclic b-chromatic number $A_b(G)$ of $G$ as a kind of a dual to acyclic chromatic number $A(G)$ of $G$.  

We continue with a similar idea to investigate the dual of some special chromatic number of a graph $G$ as the maximum number of colors in a minimal element of the mentioned strict order on all special colorings. In this contribution we concentrate on the star chromatic number $S(G)$ of a graph $G$ and introduce analogously the star b-chromatic number $S_b(G)$. 

In the next section we set the notation and give some information on the literature. In the following section we define the star b-chromatic number and present some basic results, in particular, we define the star b-chromatic number and present a natural---but more complicated than for b-chromatic number---characterization of the star b-chromatic number. In Section \ref{secdegree} we present the star degree of a vertex which leads to a natural upper bound for $S_b(G)$ presented in Section 5 and to some exact results for $S_b(G)$ for some families of graphs in Section 6. This yields further some families where $S_b(G)$ differs from other related invariants by an arbitrary big difference.  
   

\section{Preliminaries}\label{sec_Preliminaries}

In this work we consider only finite simple graphs $G=(V(G),E(G))$. We use $n_G=\vert V(G)\vert$ and $m_G=\vert E(G)\vert$. The \emph{open neighborhood} $N_G(v)$ of a vertex $v$ is the set $\{u\in V(G):uv\in E(G)\}$ and $N_G[v]=N_G(v)\cup\{v\}$ is the \emph{closed neighborhood}. The \emph{degree} $d_G(v)$ of $v\in V(G)$ is the cardinality of $N_G(v)$ and $\Delta (G)$ and $\delta (G)$ are the maximum and the minimum degree of a vertex in $G$, respectively. The clique number of $G$ is denoted by $\omega (G)$.  For $S\subseteq V(G)$ we denote by $G[S]$ the subgraph of $G$ induced on $S$. Graph $\overline{G}$ is the complement of $G$. As usual, the distance between two vertices $u$ and $v$ is denoted by $d(u,v)$ and is the minimum number of edges on a path between $u$ and $v$. The \emph{girth} $g(G)$ of a graph $G$ is the length of the shortest cycle of $G$ or infinite if such a cycle does not exists. The \emph{girth of a vertex} $v$, denoted by $g(v)$, is the length of a shortest cycle that contains $v$ or infinite if such a cycle does not exists. We use $[k]$ to denote the set $\{1,\dots,k\}$ and $[j,k]$ to denote the set $\{j,\dots,k\}$, in particular, $[1,k]=[k]$.

We call a map $c:V(G)\rightarrow [k]$ a \emph{proper vertex $k$-coloring} if $c(x)\neq c(y)$ for every edge $xy\in E(G)$. The terms "proper" and "vertex" will be omited, as we consider only proper vertex $k$-colorings. So, we call $c$ a $k$-coloring or just a coloring of $G$ in the remainder of the paper. The \emph{trivial coloring} of $G$ is the coloring with $n_G$ colors, which means that every vertex obtains a different color. The minimum number $k$ of colors, for which there exists a $k$-coloring, is called the \emph{chromatic number} $\chi(G)$ of $G$. Every $k$-coloring $c$ generates a partition of $V(G)$ into independent sets $V_i=\{u\in V(G):c(u)=i\}$, $i\in [k]$, called \emph{color classes} of $c$. We use $V_{i,j}=V_i\cup V_j$ and by $V_{i,j,\ell}=V_i\cup V_j\cup V_{\ell}$ for any $i,j,\ell\in [k]$. In particular, we use $V_{i,j,\ell}(v)$ for the component of $G[V_{i,j,\ell}]$ that contains vertex $v$. By $CN_c(v)$ we denote the set of all the colors that are present in $N_G(v)$ under coloring $c$, that is $CN_c(v)=\{c(u):u\in N_G(v)\}$. In addition we have $CN_c[v]=CN_c(v)\cup\{c(v)\}$.

\subsection{Star chromatic number}

The term star chromatic number requires careful interpretation. Namely, there are two different concepts in the literature, both with the same name: star chromatic number. First one was introduced in 1988 by Vince \cite{Vinc}. This version was later renamed to circular chromatic number by Zhu, see his survey \cite{Zhu}, and this name is still (mainly) in use today.

The second version of the star chromatic number $S(G)$ of a graph $G$ will be used from now on in this paper. It was introduced by Fertin et al. \cite{FeRR1} in 2001 as a minimum number of colors in a $k$-coloring of $G$ with the additional condition that $G[V_{i,j}]$ is a forest without a path $P_4$ for every $i,j\in [k]$. A coloring with this property is called a \emph{star coloring}. In other words, vertices of any two colors of some star coloring of $G$ induce a forest where every component is a star and $S(G)$ is the minimum number of colors in a star coloring. 

The motivation for $S(G)$ origin in the acyclic chromatic number $A(G)$ of a graph $G$, which is a minimum number of colors in a proper coloring where vertices of any two color classes induce a forest, which was introduced by Gr\"unbaum \cite{Grun}. In the same paper Gr\"unbaum already suggested that one could consider colorings where every forest induced by two color classes consists of stars. Every star coloring is also an acyclic coloring and $S(G)$ is an upper bound of more established acyclic chromatic number $A(G)$. 
 
In the seminal paper \cite{FeRR1}, see also \cite{FeRR2}, $S(G)$ and some bounds for $S(G)$ where derived for some well known families of graphs, together with some asymptotic bounds with respect to $\Delta(G)$. Ne\v set\v ril and Ossona de Mendez \cite{NeOM} improved the bound for $S(G)$ for planar graphs to $30$. Albertson et al. \cite{ACKKR} further improved $30$ to $20$ and present a planar graph with $S(G)=10$. Beside that they also improve some other bounds from \cite{FeRR1}, in particular for graphs with bounded tree-width to the best possible bound. Also a complexity question was answered in \cite{ACKKR}, where the problem of showing that $G$ has a star coloring with $3$ colors is an NP-hard problem even if $G$ is a bipartite planar graph. 

Later the investigation of star colorings followed into different directions. The most popular was the study of $S(G)$ for certain graph classes, see \cite{KaMo,KiKT} for a flavor. Star chromatic index was introduced by Dvo\v r\'ak et al. \cite{DvMS} and later its list version by Lu\v zar et al. \cite{LuMS}. Some complexity and algorithmic issues were considered by Shalu and Antony \cite{ShAn} and  Gebremedhin et al. \cite{GTMP}, respectively.

\subsection{b-chromatic number}

We start with the original definition of the b-chromatic number from \cite{IrMa}. By $\mathcal{F}(G)$ we denote the set of all colorings of a graph $G$ and choose a $k$-coloring $c\in \mathcal{F}(G)$. Let $v\in V(G)$ with $c(v)=i$. If $CN_c[v]=[k]$, that is, if all the $k$ colors are present in the closed neighborhood of $v$, then $v$ is a b-\emph{vertex} (of color $i$). If at least one color, say $j$, is missing in $N_G[v]$, then $v$ is not a b-vertex and we can recolor $v$ with $j$ to obtain a slightly different coloring. Moreover, if there exists no b-vertex of color $i$, then we can recolor every vertex $v$ of color $i$ with some color not present in $CN_c[v]$, say $j_v$. Therefore $c_i:V(G)\rightarrow [k]-\{i\}$ defined by 
\begin{equation*}
c_i(v)=\left\{ 
\begin{array}{ccc}
c(v) & : & c(v)\neq i \\
j_v & : & c(v)=i 
\end{array}
\right. 
\end{equation*}
is a $(k-1)$-coloring of $G$. This procedure is called a \emph{recoloring step}. An iterative performing of recoloring steps from a trivial coloring of $G$ while it is possible is called the \emph{recoloring algorithm}. (Actually one could start with any coloring from $\mathcal{F}(G)$.) It was mentioned in \cite{AnCP}, see subsection 2.2, that the recoloring algorithm has a polynomial time complexity and can be therefore considered as a heuristic approach to the chromatic number $\chi(G)$. 

Now we define the relation $\triangleleft$ on $\mathcal{F}(G)\times \mathcal{F}(G)$. If $c'\in \mathcal{F}(G)$ can be obtained from $c\in \mathcal{F}(G)$ by a recoloring step of one fixed color class of $c$, then $c'$ is in relation $\triangleleft$ with $c$, that is $c'\triangleleft c$.
Relation $\triangleleft$ is asymmetric because $c'$ has less colors than $c$. Let $\prec$ be the transitive closure of $\triangleleft$, which is a strict partial order on $\mathcal{F}(G)$. For a finite graph $G$ we have finitely many different colorings of $G$ and $\prec$ has some minimal elements. Now, the chromatic number $\chi (G)$ is the minimum number of colors of a minimal element of $\prec$ and a dual of $\chi(G)$ is the maximum number $\varphi(G)$ of colors of a minimal element of $\prec$.  We call $\varphi (G)$ the \emph{b-chromatic number} of $G$. 

A minimal element of $\prec$ is called a \emph{b-coloring}. We cannot perform a recoloring step on a minimal element of $\prec$, which means that every color class of any minimal element of $\prec$ has a b-vertex. Hence, the alternative definition, as already mentioned in \cite{IrMa}, is that the b-chromatic number is the maximum number of colors in a b-coloring of $G$. This definition was later used in almost all publications on the b-chromatic number and related topics. 

Every b-coloring with $\varphi(G)$ colors needs at least $\varphi(G)$ vertices of high enough degree to be b-vertices (at least one for each color). The degree of a b-vertex in such a coloring must be $d_G(v)\geq\varphi(G)-1$. Let vertices $v_1,\dots,v_{n_G}$ of $G$ be ordered by degrees $d_G(v_1)\geq \cdots\geq d_G(v_{n_G})$. The $m$-\emph{degree} $m(G)$ is defined as 
$$m(G)=\max \{i: i-1\leq d_G(v_i)\}.$$
With $m$-degree one obtain a natural upper bound for the b-chromatic number: $\varphi(G)\leq m(G)$, as already shown in \cite{IrMa}. 

Determining of the b-chromatic number of a graph is an NP-hard problem \cite{IrMa}. The problem remains NP-hard for connected bipartite graphs see Kratochv{\'{\i}}l et al. \cite{krtuvo-02}. In contrast, the exact result for every tree $T$ was given already in \cite{IrMa} and is very close to $m(G)$. This was later generalized for cactus graphs \cite{CaLSMaSi}, for outerplanar graphs \cite{MaSi12}, and for graphs with large enough girth \cite{CaFaSi,CaLiSi,koza-15}. We recommend survey \cite{JaPe} for more information about the b-chromatic number and related concepts.



\section{Definition and some basic results}\label{sec_Definition}

In this section we define the star b-chromatic number $S_b(G)$ in an analogue fashion as the acyclic b-chromatic number $A_b(G)$ was defined in \cite{AnCP}. As explained for $A_b(G)$ and $\varphi(G)$ in \cite{AnCP} we cannot observe the intersection of all b-colorings and all star colorings, because there are graphs where this intersection is empty. The smallest example for this is $C_4$, but it also holds for the whole family of complete bipartite graphs $K_{p,q}$ where ${\rm min}\{p,q\}\geq 2$. Similar one can not define $S_b(G)$ from $A_b(G)$ because the only acyclic b-coloring of $P_4$ is not a star coloring. We can avoid this problems if we follow the original definition of $\varphi(G)$, only that we need to restrict from all colorings of $G$ to all star colorings of $G$.

By $\mathcal{SF}(G)$ we denote the set of all star colorings of a graph $G$. A recoloring step for $c\in\mathcal{SF}(G)$ that produces coloring $c'$ is a \emph{star recoloring step} if $c'\in \mathcal{SF}(G)$. So, in a star recoloring step we can reduce some color of a star coloring $c$ only when a new coloring $c'$ is also a star coloring.
The process of applying a star recoloring step to a trivial coloring until no further steps are possible is called the \emph{star recoloring algorithm}. Note that the star recoloring algorithm has polynomial time complexity. For this we need, in addition to the recoloring algorithm, to check whether the new coloring is still star coloring after every star recoloring step. There are at most $O(n_G^2)$ pairs of different colors and one can do this in polynomial time. Therefore, one can observe the star recoloring algorithm as a heuristic algorithm for an approximation value of the star chromatic number. 

Now we define a relation $\triangleleft_s\subseteq\mathcal{SF}(G)\times \mathcal{SF}(G)$. For $c',c\in \mathcal{SF}(G)$ we have $c'\triangleleft_s c$ when $c'$ is obtained from $c$ by a star recoloring step. Also $\triangleleft_s$ is asymmetric because $c'$ has less colors than $c$. By $\prec_s$ we denote its transitive closure. Therefore $\prec_s$ is a strict partial order of $\mathcal{SF}(G)$.  Clearly, $\mathcal{SF}(G)$ is finite as $G$ is finite and $\prec_s$ contains minimal elements. The following proposition shows that the trivial coloring $t$ is the greatest element of $\prec_s$ (sometimes also called the maximum element). 

\begin{proposition}\label{first}
 Let $G$ be a graph and $t\in\mathcal{SF}(G)$ be a trivial coloring. If $c\in\mathcal{SF}(G)$, then there exists a chain 
 $$c\triangleleft_s c_1\triangleleft_s c_2\triangleleft_s\cdots\triangleleft_s c_{\ell-1}\triangleleft_s t.$$
 \end{proposition}
 
\noindent {\textbf{Proof.}}
 Let $c\in \mathcal{SF}(G)$ be an arbitrary star $k$-coloring. Let $\ell=n_G-k$, $c_{\ell}=t$ and $c_0=c$. We may assume that the first $k$ colors from $t$ are the $k$ colors from $c$. We can order $v_1,\dots,v_k,v_{k+1},\dots,v_{n_G}$ the vertices of $G$ such that $c_0(v_i)=i=c_{\ell }(v_i)$ for every $i\in [k]$, and the rest of the ordering is arbitrary. We define the coloring $c_i$ from $c_{i-1}$ for every $i\in [\ell]$ by 
 \begin{equation*}
c_i(v)=\left\{ 
\begin{array}{lcl}
c_{i-1}(v) & \text{if} & v\neq v_{k+i}, \\
k+i & \text{if} & v=v_{k+i}. 
\end{array}
\right.
\end{equation*}
From the other perspective, if we reverse the order of colorings, then we obtain $c$ from $t$ by recoloring every vertex $v_{n_G},\dots,v_{k+1}$ exactly once. Clearly, $c_i\in\mathcal{SF}(G)$ for every $i\in [\ell]$ because $c\in\mathcal{SF}(G)$. The color $c_{i-1}(v_{k+i})$ is not in the closed neighborhood of $v_{k+i}$ in coloring $c_i$ for every $i\in [\ell]$. Therefore, one can perform a star recoloring step from $c_i$ to $c_{i-1}$ for every $i\in [\ell]$.  Hence, $c_{i-1}\triangleleft_s c_i$ follows for every $i\in [\ell]$. \qed\bigskip

With this the following definition is justified. The \emph{star b-chromatic number} $S_b(G)$ is the maximum number of colors in a minimal element of $\prec_s$. In other words,
$$S_b(G)=\max\{|c|:c\in\mathcal{SF}(G) \text{ is a minimal element of }\prec_s\}.$$
The star b-chromatic number of a graph $G$ describes the worst case to appear while using the star recoloring algorithm to estimate $S(G)$. A star coloring of $G$ with $S_b(G)$ colors that arise from a minimal element of $\prec_s$ is called an $S_b(G)$-\emph{coloring}. We have the following inequality chain

\begin{equation}\label{chain}
\omega(G)\leq\chi(G)\leq S(G)\leq S_b(G)\leq n_G,
\end{equation}
where the fact that every star b-coloring is also a star coloring implies $S(G)\leq S_b(G)$. 
Next we characterize the graphs for which $S_b(G)=n_G$. 

\begin{proposition}\label{order}
We have $S_b(G)=n_G$ if and only if $G\cong K_{n_G}$. 
\end{proposition}

\noindent {\textbf{Proof.}}
If $G\cong K_n$, then $n=\chi(K_n)\leq S_b(K_n)\leq n$ by (\ref{chain}). Conversely, if $G\ncong K_n$, then there exist different and nonadjacent $u,w\in V(G)$. For a trivial coloring $t$ of $G$,
 \begin{equation*}
c(v)=\left\{ 
\begin{array}{lccl}
t(v) &:& \text{if} & v\neq u, \\
t(w) &:& \text{if} & v=u, 
\end{array}
\right. 
\end{equation*}
is a coloring obtained by a star recoloring step. Hence $t$ is not a minimal element of $\prec_s$ and $S_b(G)<n_G$.
\qed\bigskip

Next we describe the vertices that cannot be recolored in a star coloring $c=(V_1,\dots,V_k)$ of a graph $G$. A color $j\in [k]$ is \emph{blocked} for $v\in V_i$ (with respect to $c$) if 
 \begin{equation*}
c'(u)=\left\{ 
\begin{array}{lccl}
c(u) &:& \text{if} & u\neq v, \\
j &:& \text{if} & u=v, 
\end{array}
\right. 
\end{equation*} 
is not a star coloring of $G$. Beside that, also the color $i$ is blocked for $v\in V_i$. Clearly, all the colors from $CN_c[v]$ are blocked for $v$. Further, a color $\ell\in[k]-CN_c[v]$ is blocked for $v$ if there exists a color $j\in CN_c(v)$, such that $G[V_{j,\ell}\cup\{v\}]$ contains a $P_4$, because this $P_4$ is bi-colored if we recolor $v$ by $\ell$. If every color from $[k]$ is blocked for $v$, then $v$ cannot be recolored and with this there exists no star recoloring step for $V_i$. Therefore we call $v\in V_i$ a \emph{strong star b-vertex} for $c$ if every color from $[k]-\{i\}$ is blocked for $v$. Obviously, every b-vertex is a strong star b-vertex as well. Moreover, this is a generalization of b-vertices, because no star recoloring step is possible if every color class contains a strong star b-vertex.  

As already suggested by "strong" in strong star b-vertex, this does not generalize b-vertices to star b-vertices---our goal---in all the cases. The simplest example to observe this is a five-cycle with consecutive colors $1-4-2-3-4$. It is straightforward to check that vertices of colors $1,2$ and $3$ are strong star b-vertices of respective color, however there is no strong star b-vertex of color $4$. Nevertheless, we are not able to execute a star recoloring step for color $4$, because $3$ ($2$, respectively) is the only color that is not blocked for the first (second, respectively) vertex of color $4$ and such a recoloring yields a bi-colored $P_4$. 

To generalize this situation, notice that we have a $P_4$ with colors $\alpha \beta \gamma \alpha $ where first vertex of color $\alpha $ can be recolored only by $\gamma $ and the second vertex of color $\alpha $ can be recolored only with color $\beta $. Clearly, such a path must be induced. The situation can be more general and a color $j\in [k]$ that is not blocked for $v\in V_i$ is called an \emph{available color} for $v$. Further, $u,v\in V_i$ are in relation $\sim_i$ if there exist an induced $P_4$ that starts in $u$ and ends in $v$. Relation $\sim_i$ is symmetric but not reflexive and transitive. Therefore we consider reflexive-transitive closure $\sim_i^*$ of $\sim_i$, which is an equivalence relation. We call the equivalence class of a vertex $v\in V_i$ with respect to $\sim_i^*$ the $P_4$-\emph{system} of $v$ and denote it shortly by $P_4(v)$. An available color $j\in [k]$ of $v$ is \emph{blocked} for $P_4(v)$ if any recoloring of $V_i$ where $v$ is recolored by $j$ yields a bi-colored $P_4$ between two vertices of $P_4(v)$. Often $P_4(v)$ is a single vertex and in such a case $P_4(v)$ is not blocked for any available color of $v$. On the other hand, for the mentioned coloring of $C_5$, let $u$ and $v$ be the first and second vertex, respectively, of color 4. Clearly, $P_4(u)=\{u,v\}=P_4(v)$ and only available color for $u$ is $3$ and only available color for $v$ is $2$. Now, $3$ is blocked for $P_4(u)$ and symmetrically $2$ is blocked for $P_4(v)$.

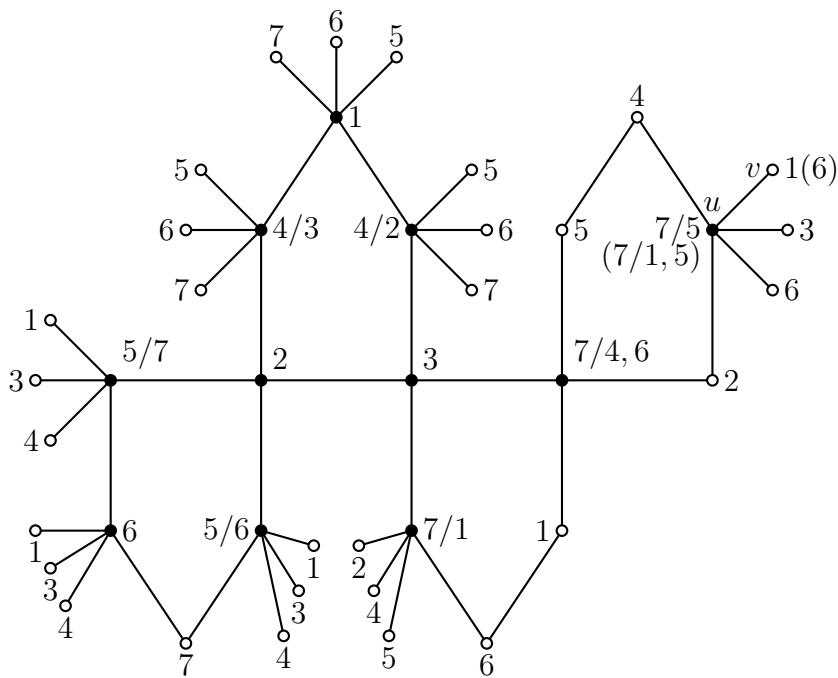
\begin{figure}[ht!]
\begin{center}
\begin{tikzpicture}[xscale=1, yscale=1,style=thick,x=1cm,y=1cm]
\def\vr{2pt} 

\path (0,0) coordinate (a);
\path (2,0) coordinate (b);
\path (4,0) coordinate (c);
\path (6,0) coordinate (d);
\path (8,0) coordinate (e);
\path (0,-2) coordinate (f);
\path (2,-2) coordinate (g);
\path (4,-2) coordinate (h);
\path (6,-2) coordinate (i);
\path (2,2) coordinate (j);
\path (4,2) coordinate (k);
\path (6,2) coordinate (l);
\path (8,2) coordinate (m);
\path (1,-3.5) coordinate (n);
\path (5,-3.5) coordinate (p);
\path (3,3.5) coordinate (q);
\path (7,3.5) coordinate (r);

\path (-1,0) coordinate (a1);
\path (-0.8,0.8) coordinate (a2);
\path (-0.8,-0.8) coordinate (a3);

\path (9,2) coordinate (m1);
\path (8.8,2.8) coordinate (m2);
\path (8.8,1.2) coordinate (m3);

\path (5,2) coordinate (k1);
\path (4.8,2.8) coordinate (k2);
\path (4.8,1.2) coordinate (k3);

\path (1,2) coordinate (j1);
\path (1.2,2.8) coordinate (j2);
\path (1.2,1.2) coordinate (j3);

\path (3,4.5) coordinate (q1);
\path (3.8,4.3) coordinate (q2);
\path (2.2,4.3) coordinate (q3);

\path (-1,-2) coordinate (f1);
\path (-0.8,-2.5) coordinate (f2);
\path (-0.6,-3) coordinate (f3);

\path (2.7,-2.2) coordinate (g1);
\path (2.5,-2.8) coordinate (g2);
\path (2.3,-3.4) coordinate (g3);

\path (3.3,-2.2) coordinate (h1);
\path (3.5,-2.8) coordinate (h2);
\path (3.7,-3.4) coordinate (h3);


\draw (a1) -- (a) -- (a2); 
\draw (m1) -- (m) -- (m2); 
\draw (j1) -- (j) -- (j2); 
\draw (q1) -- (q) -- (q2); 
\draw (f1) -- (f) -- (f2); 
\draw (g1) -- (g) -- (g2); 
\draw (h1) -- (h) -- (h2); 
\draw (k1) -- (k) -- (k2); 
\draw (a3) -- (a) -- (f) -- (f3); 
\draw (f) -- (n) -- (g) -- (g3);
\draw (a) -- (b) -- (j) -- (j3); 
\draw (g) -- (b) -- (c) -- (h) -- (h3);  
\draw (h) -- (p) -- (i) -- (d); 
\draw (m3) -- (m) -- (e) -- (d) -- (l) -- (r) -- (m); 
\draw (d) -- (c) -- (k) -- (k3);
\draw (k) -- (q) -- (q3);
\draw (j) -- (q);

\draw (a) [fill=black] circle (\vr);
\draw (b) [fill=black] circle (\vr);
\draw (c) [fill=black] circle (\vr);
\draw (d) [fill=black] circle (\vr);
\draw (e) [fill=white] circle (\vr);
\draw (f) [fill=black] circle (\vr);
\draw (g) [fill=black] circle (\vr);
\draw (h) [fill=black] circle (\vr);
\draw (i) [fill=white] circle (\vr);
\draw (j) [fill=black] circle (\vr);
\draw (k) [fill=black] circle (\vr);
\draw (l) [fill=white] circle (\vr);
\draw (m) [fill=black] circle (\vr);
\draw (n) [fill=white] circle (\vr);
\draw (p) [fill=white] circle (\vr);
\draw (q) [fill=black] circle (\vr);
\draw (r) [fill=white] circle (\vr);
\draw (a1) [fill=white] circle (\vr);
\draw (a2) [fill=white] circle (\vr);
\draw (a3) [fill=white] circle (\vr);
\draw (f1) [fill=white] circle (\vr);
\draw (f2) [fill=white] circle (\vr);
\draw (f3) [fill=white] circle (\vr);
\draw (g1) [fill=white] circle (\vr);
\draw (g2) [fill=white] circle (\vr);
\draw (g3) [fill=white] circle (\vr);
\draw (h1) [fill=white] circle (\vr);
\draw (h2) [fill=white] circle (\vr);
\draw (h3) [fill=white] circle (\vr);
\draw (j1) [fill=white] circle (\vr);
\draw (j2) [fill=white] circle (\vr);
\draw (j3) [fill=white] circle (\vr);
\draw (q1) [fill=white] circle (\vr);
\draw (q2) [fill=white] circle (\vr);
\draw (q3) [fill=white] circle (\vr);
\draw (k1) [fill=white] circle (\vr);
\draw (k2) [fill=white] circle (\vr);
\draw (k3) [fill=white] circle (\vr);
\draw (m1) [fill=white] circle (\vr);
\draw (m2) [fill=white] circle (\vr);
\draw (m3) [fill=white] circle (\vr);

\draw[anchor = south west] (a) node {$5/7$};
\draw[anchor = south west] (b) node {$2$};
\draw[anchor = south west] (c) node {$3$};
\draw[anchor = south west] (d) node {$7/4,6$};
\draw[anchor = west] (e) node {$2$};
\draw[anchor = west] (f) node {$6$};
\draw[anchor = east] (g) node {$5/6$};
\draw[anchor = west] (h) node {$7/1$};
\draw[anchor = east] (i) node {$1$};
\draw[anchor = west] (j) node {$4/3$};
\draw[anchor = east] (k) node {$4/2$};
\draw[anchor = west] (l) node {$5$};
\draw[anchor = east] (m) node {$7/5$};
\draw[anchor = north east] (m) node {$(7/1,5)$};
\draw[anchor = south] (m) node[yshift=1mm] {$u$};
\draw[anchor = north] (n) node {$7$};
\draw[anchor = north] (p) node {$6$};
\draw[anchor = west] (q) node {$1$};
\draw[anchor = south] (r) node {$4$};
\draw[anchor = east] (a1) node {$3$};
\draw[anchor = east] (a2) node {$1$};
\draw[anchor = east] (a3) node {$4$};
\draw[anchor = north] (f1) node {$1$};
\draw[anchor = north] (f2) node {$3$};
\draw[anchor = north] (f3) node {$4$};
\draw[anchor = north] (g1) node {$1$};
\draw[anchor = north] (g2) node {$3$};
\draw[anchor = north] (g3) node {$4$};
\draw[anchor = north] (h1) node {$2$};
\draw[anchor = north] (h2) node {$4$};
\draw[anchor = north] (h3) node {$5$};
\draw[anchor = east] (j1) node {$6$};
\draw[anchor = east] (j2) node {$5$};
\draw[anchor = east] (j3) node {$7$};
\draw[anchor = south] (q1) node {$6$};
\draw[anchor = south] (q2) node {$5$};
\draw[anchor = south] (q3) node {$7$};
\draw[anchor = west] (k1) node {$6$};
\draw[anchor = west] (k2) node {$5$};
\draw[anchor = west] (k3) node {$7$};
\draw[anchor = west] (m1) node {$3$};
\draw[anchor = west] (m2) node {$1(6)$};
\draw[anchor = east] (m2) node {$v$};
\draw[anchor = west] (m3) node {$6$};

\end{tikzpicture}
\end{center}
\caption{
Two star colorings of a graph. In the first coloring the black vertices have all the other colors blocked, meaning that every color class has at least one star b-vertex. For the vertices that are not strong star b-vertices, the available colors are indicated. The second coloring differs from the first only in $v$, which is indicated in brackets. This also changes the available colors for vertex $u$, allowing us to apply a star recoloring step for color $7$.}
\label{blocked}
\end{figure}

For further example consider Figure \ref{blocked} that contains a graph $G$ together with two colorings-the second coloring is almost the same as the first, only vertex $v$ has color $6$ instead of $1$ (denoted by $1(6)$ on figure). By $i/j$ we denote a color $i$ of some vertex where color $j$ is available for that vertex. For the first coloring notice that all colors are blocked for black vertices. The black vertices of colors $1,2,3$ and $6$ are unique and are strong star b-vertices. There are two black vertices of color $4$ and of color $5$ together with their available colors. Notice that each pair form a $P_4$-system for $4$ and $5$, respectively, and that the available colors are blocked in each case. Further, there are three black vertices of color $7$ together with their available colors, which are also blocked in every case. The second coloring is almost the same as the first one, only vertex $v$ has color $6$ instead of $1$. With this $1$ became an available color for $u$ and $P_4$-system of black vertices of color $7$ is not blocked anymore. It is easy to run a star recoloring step for color $7$ of the second coloring.   

Now we can define a star b-vertex that is a star analogue to b-vertex and will further allow to describe all minimal elements of $\prec_s$. Unfortunately, it is not so elegant as the result for $\varphi(G)$, see \cite{IrMa}, and is more similar to acyclic b-vertices and the description of $A_b(G)$, see \cite{AnCP}. 

\begin{definition}
Let $G$ be a graph with a star coloring $c:V(G)\rightarrow [k]$. A vertex $v\in V_i$, $i\in [k]$, is a star b-vertex if $P_4(v)$ is blocked for any available color for $v$ with respect to $c$. 
\end{definition}

Every b-vertex $v$ is also a star b-vertex, since in this case  $[k]-CN_c[v]=\emptyset$ and there are no available colors for $v$. Similarly there are no available colors for $v$ when $v$ is a (strong) star b-vertex. Moreover, the minimal elements of $\prec_s$ of a graph $G$ must contain a star b-vertex of every color as shown next. This further leads to the alternative description of $S_b(G)$.

\begin{theorem}\label{conditionABcoloring}
A star $k$-coloring $c$ is a minimal element of $\prec_s$ if and only if every color class $V_i$, $i\in[k]$, contains a star b-vertex.
\end{theorem}

\noindent {\textbf{Proof.}}
Let $G$ be a graph and let a star $k$-coloring $c$ be a minimal element of $\prec_s$ of $G$. Thus a star recoloring step cannot be performed for $c$. If $V_i$, $i\in [k]$, has a (strong star) b-vertex $x$, then $x$ is also a star b-vertex and we are done. So, assume that $V_i$ has no strong star b-vertex. Let $v\in V_i$ be a vertex for which the star recoloring step fails. We may choose $v$ in such a way that the minimum number $t$ of vertices from $V_i$ must be recolored beside $v$, such that we get a bi-colored $P_4$. Let $\ell$ be any available color for $v$, which means that the coloring 
 \begin{equation*}
c'(u)=\left\{ 
\begin{array}{lccl}
c(u) &:& \text{if} & u\neq v, \\
\ell &:& \text{if} & u=v, 
\end{array}
\right. 
\end{equation*}
is a star coloring of $G$. Hence, we get a bi-colored $P_4$ after also some others vertices of $V_i$ are recolored. By induction on $t$ we proceed to show that we get a system $P_4(v)$. Let first $t=1$ and let $u\in V_i$ be the other vertex that needs to be recolored. Clearly, $uv\notin E(G)$ and $2\leq d(u,v)\leq 3$. Suppose first that there is no induced $P_4$ between $u$ and $v$. The only $P_4$ that contains both $u$ and $v$ is of the form $vxuy$ or $xvyu$. But then $c(x)=c(y)$ and both $u$ and $v$ must be recolored by the same color to obtain bi-colored $P_4$, a contradiction since $c(x)\neq c(y)$ as $c$ is a star coloring. This means that there exists an induced $P_4=vxyu$ and $c(y)=\ell $ is the only available color for $v$ and $c(x)$ the only available color for $u$ since $t=1$. In other words, color $\ell$ is blocked by $P_4(v)$. 

Let now $t>1$ and let $v=v_0,v_1,\dots,v_t\in V_i$ be the minimum number of vertices that need to be recolored to obtain a bi-colored $P_4$. By the same reasons as above we get that $d(v_0,v_1)=3$ and that there exists an induced path $v_0xyv_1$ with $c(y)=\ell$ and by $c'$ there are some other available colors for $v_1$ beside the color $c(x)$. By induction hypothesis there exists a system $P_4(v_1)$ that yields a bi-colored $P_4$ when $c(v_1)\neq c(x)$ or a bi-colored $v_0xyv_1$ when $c(v_1)=c(x)$. So, $V_i$ has a $v_1$ star b-vertex and we are done with this direction.

Conversely, if every color class $V_i$, $i\in [k]$, of $k$-coloring $c$ contains a star b-vertex, say $v_i$, then $v_i$ cannot be recolored by any color for every $i\in [k]$ and $c$ is a minimal element of $\prec_s$. 
\qed\bigskip

\begin{corollary}\label{cor}
The star b-chromatic number $S_b(G)$ of a graph $G$ is the largest integer $k$, such that there exists a star $k$-coloring, where every color class $V_i$, $i\in [k]$, contains a star b-vertex.
\end{corollary}


\section{Star degree of a vertex}\label{secdegree} 

In this section we define the star degree of a vertex which yields the star m-degree of a graph analogous to m-degree for b-chromatic number and acyclic m-degree for acyclic b-chromatic number, for the later see \cite{AnCP}. The idea is in a (star, acyclic) coloring of a graph $G$ that represents the maximum number of colors that can be blocked for a vertex $v$. For the b-chromatic number is this simply the degree of a vertex, hence the name. It is more complicated for acyclic degree of a vertex, but still relatively easy to describe, not always easy to compute it, see \cite{AnCP}. As we will see in this section, the star degree seems to be much more complicated to describe in general situation, but easy enough to deal with for vertices that do not belong to a three or a four cycle. 

The \emph{star degree} $d_G^s(v)$ of a vertex $v$ of a graph $G$ is the maximum number of colors that can be blocked for $v$ in a star coloring $c$ of $G$. All the colors that appear in $N_G(v)$ are clearly blocked. The second possibility with respect to $v$ represent paths $P_4=xvyu$ or $P_4=vxuy$ where $c(x)=c(y)$ and the color $c(u)\neq c(v)$ is blocked for $v$. Such a path $P_4$ is called \emph{blocking} with respect to $v$. Notice that in every system $P_4(v)$ there exists at least one path $P_4$ as already described for every blocking color and we do not need to take system $P_4(v)$ into consideration here as we will only count the maximum number of blocked colors for $d_G^s(v)$.

Blocking paths $P_4$ with respect to $v$ will be further separated by the distances of vertices of $P_4$ from $v$. By this we mean that blocking path $abcd$, $a,b,c,d\in [0,3]$,  with respect to $v$ is a path $P_4=pqrs$ where $d_G(v,p)=a$, $d_G(v,q)=b$, $d_G(v,r)=c$ and $d_G(v,s)=d$. Clearly, exactly one of $p,q,r,s$ equals to $v$. By the symmetry of $P_4$ we distinguish blocking paths $0123$, $0122$, $0121$, $0112$, $0111$, $1012$ and $1011$. Moreover, blocking paths $0111$ and $1011$ contain only the neighbors of $v$ and these colors are blocked. Hence, we will ignore them. Notice also that a $0121$ blocking path induces a $C_4$ and can be therefore considered as $1012$ blocking path. Finally, in a $0112$ blocking path $vxuy$ has $c(x)=c(y)$ and $c(x)$ and $c(u)$ are blocked for $v$ since $x$ and $u$ are neighbors of $v$. Thus, such a blocking path does not bring any new color that is blocked and we can ignore them as well. So, we are left with $0123$, $0122$ and $1012$ blocking paths. \\

\begin{figure}[!ht]
	\begin{center}
	\begin{tikzpicture}[xscale=0.5, yscale=0.5]
		\tikzstyle{rn}=[circle,fill=white,draw, inner sep=0pt, minimum size=5pt]
		
		\tikzstyle{bn}=[circle,fill=black,draw, inner sep=0pt, minimum size=5pt]
		\tikzstyle{every node}=[font=\footnotesize]
		
		\node (a)[rn] at (0 cm,0 cm){};
		\node (b)[rn] at (1 cm,0 cm){};
		\node (c)[rn] at (2 cm,0 cm){};
		\node (d)[rn] at (3 cm,0 cm){};
		
		\draw (a)--(b);
		\draw (b)--(c);
		\draw (c)--(d);
		
		\draw[anchor = north] (a) node {$0$};
		\draw[anchor = south] (a) node {$v$};
		\draw[anchor = north] (b) node {$1$};
		\draw[anchor = north] (c) node {$2$};
		\draw[anchor = north] (d) node {$3$};
		
		\node (a1)[rn] at (5 cm,0 cm){};
		\node (b1)[rn] at (6 cm,0 cm){};
		\node (c1)[rn] at (7 cm,0 cm){};
		\node (d1)[rn] at (8 cm,0 cm){};
		\node (e1)[rn] at (6.5 cm,1 cm){};
		
		\draw (a1)--(b1);
		\draw (b1)--(c1);
		\draw (c1)--(d1);
		\draw (a1)--(e1);
		\draw (e1)--(d1);
		
		\draw[anchor = north] (a1) node {$0$};
		\draw[anchor = south] (a1) node {$v$};
		\draw[anchor = north] (b1) node {$1$};
		\draw[anchor = north] (c1) node {$2$};
		\draw[anchor = north] (d1) node {$2$};
		
		\node (a2)[rn] at (10 cm,0 cm){};
		\node (b2)[rn] at (11 cm,0 cm){};
		\node (c2)[rn] at (12 cm,0 cm){};
		\node (d2)[rn] at (13 cm,0 cm){};
		
		\draw (a2)--(b2);
		\draw (b2)--(c2);
		\draw (c2)--(d2);
		
		\draw[anchor = north] (a2) node {$1$};
		\draw[anchor = north] (b2) node {$0$};
		\draw[anchor = south] (b2) node {$v$};
		\draw[anchor = north] (c2) node {$1$};
		\draw[anchor = north] (d2) node {$2$};
	\end{tikzpicture}
\end{center}
\caption{Blocking paths $0123$, $0122$ and $1012$, shown from left to right.}
\end{figure}
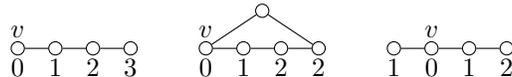

Our goal is to count the maximum number of blocked colors for $v\in V(G)$. Beside the colors from $N_G(v)$, every blocking $P_4$ with respect to $v$ brings at most one additional blocked color. Unfortunately, there are some intersections in the mentioned three cases of blocking paths and not every blocking $P_4$ yields a new blocked color. For this observe, for instance, two $0123$ blocking paths $vxuy$ and $vx_1u_1y$ where $x\neq x_1$. We have either $c(x)=c(y)=c(x_1)$ and $c(x),c(u),c(u_1)$ are blocked or $c(y)=c(x)\neq c(x_1)$ (up to the symmetry) and $c(x),c(u),c(x_1)$ are blocked. In both cases at most three colors are blocked, but not all four $c(x),c(x_1),c(u),c(u_1)$, which is maximal possible for two blocking paths.

Another situation that decrees the number of blocked colors are edges in $G[N_G(v)]$. The maximum number of colors in $N_G(v)$ is $d_G(v)$. However, sometimes it is more convenient to have less colors in $N_G(v)$ to achieve more blocked colors with blocking paths as we will see. 

All in all, it seems that it is not easy to get $d_G^s(v)$ in general. However, we can compute $d_G^s(v)$ in some special cases as we will see in the rest of this section. We start with additional condition that $g(v)\geq 7$. For this we need some further notations about vertices that form blocking paths. Denote with $N_2(v)$ the set of all vertices $u$ for which $d_G(v,u)=2$ and with $N_3(v)$ the set of all vertices $u$ for which $d_G(v,u)=3$. The set $X(v)\subseteq N_2(v)$ contains all vertices at distance two to $v$ which have no neighbor in $N_3(v)$. The set $Y(v)\subseteq N_2(v)$ contains all vertices at distance two to $v$ which have a neighbor in $N_3(v)$. A set $A_2(v)\subseteq N_G(v)$ contains the minimum number of neighbors of $v$, such that every vertex of $X(v)$ has a neighbor in $A_2(v)$. The minimality of $A_2(v)$ means that for every $u\in A_2(v)$   there exists $x_u\in X(v)\cap N_G(u)$ such that $N_G(x_u)\cap N_G(v)=\{u\}$. Therefore we call every vertex $x\in X(v)$ with $N_G(x)\cap N_G(v)=\{u\}$ a \emph{privat neighbor} of $u$ in $X(v)$. 

\begin{figure}[ht!]
\begin{center}
\begin{tikzpicture}[xscale=1, yscale=1,style=thick,x=1cm,y=1cm, rotate=90]
\def\vr{2pt} 

\path (2,2) coordinate (v);
\path (0,2) coordinate (b);
\path (0,3.5) coordinate (c);
\path (0,5) coordinate (d);

\path (4,2) coordinate (b1);
\path (5,2) coordinate (b2);
\path (6,2) coordinate (b3);

\path (4,3.5) coordinate (c1);
\path (4,4.5) coordinate (c2);
\path (4,5.5) coordinate (c3);

\path (5,3.5) coordinate (c4);

\path (1,0) coordinate (d1);
\path (-1,0) coordinate (d2);
\path (0,-1) coordinate (d3);
\path (0,-2) coordinate (d4);

\path (2,-1) coordinate (d5);
\path (-1,-3) coordinate (d6);
\path (1,-3) coordinate (d7);

\path (3,0) coordinate (e1);
\path (5,0) coordinate (e2);
\path (4,-1) coordinate (e3);
\path (4,-2) coordinate (e4);

\path (2,3.5) coordinate (f1);
\path (2,4.5) coordinate (f2);
\path (2,5.5) coordinate (f3);
\path (1,3.5) coordinate (f4);

\draw (d) -- (c) -- (v) -- (c1) -- (c2) -- (c3);
\draw (c1) -- (c4);
\draw (b) -- (v) -- (b1) -- (b2) -- (b3);

\draw (d1) -- (v) -- (d2);
\draw (e1) -- (v) -- (e2);
\draw (d1) -- (d3) -- (d2);
\draw (e1) -- (e3) -- (e2);
\draw (d1) -- (d4) -- (d2);
\draw (e1) -- (e4) -- (e2);
\draw (d3) -- (d5);
\draw (d6) -- (d4) -- (d7);

\draw (b) -- (v) -- (b1) -- (b2) -- (b3);

\draw (v) -- (f1) -- (f2) -- (f3);
\draw (f1) -- (f4);

\draw (b) [fill=white] circle (\vr);
\draw (c) [fill=white] circle (\vr);
\draw (d) [fill=white] circle (\vr);
\draw (v) [fill=white] circle (\vr);
\draw (b1) [fill=white] circle (\vr);
\draw (b2) [fill=white] circle (\vr);
\draw (b3) [fill=white] circle (\vr);
\draw (c1) [fill=white] circle (\vr);
\draw (c2) [fill=white] circle (\vr);
\draw (c3) [fill=white] circle (\vr);
\draw (c4) [fill=white] circle (\vr);

\draw (d1) [fill=white] circle (\vr);
\draw (d2) [fill=white] circle (\vr);
\draw (d3) [fill=white] circle (\vr);
\draw (d4) [fill=white] circle (\vr);

\draw (d5) [fill=white] circle (\vr);
\draw (d6) [fill=white] circle (\vr);
\draw (d7) [fill=white] circle (\vr);

\draw (e1) [fill=white] circle (\vr);
\draw (e2) [fill=white] circle (\vr);
\draw (e3) [fill=white] circle (\vr);
\draw (e4) [fill=white] circle (\vr);

\draw (f1) [fill=white] circle (\vr);
\draw (f2) [fill=white] circle (\vr);
\draw (f3) [fill=white] circle (\vr);
\draw (f4) [fill=white] circle (\vr);

\draw[anchor = west] (v) node {$\, \, \, v$};
\draw[anchor = north] (b) node {$a_1^1$};
\draw[anchor = north] (c) node {$a_2^2$};
\draw[anchor = north] (c1) node {$a_4^2 \,$};
\draw[anchor = east] (b1) node {$a_5^3$};
\draw[anchor = north] (e1) node {$a_7^2$};
\draw[anchor = south] (e2) node {$a_6^1$};
\draw[anchor = south] (d1) node {$a_8^3$};
\draw[anchor = north] (d2) node {$a_9^1$};

\draw[anchor = south] (f1) node {$a_3^2$};
\draw[anchor = east] (f4) node {$b_2$};
\draw[anchor = south] (f2) node {$b_3$};
\draw[anchor = south] (f3) node {$c_1$};

\draw[anchor = north] (d) node {$b_1$};
\draw[anchor = north] (c2) node {$b_4$};
\draw[anchor = north] (c3) node {$c_2$};
\draw[anchor = south] (c4) node {$b_5$};

\draw[anchor = east] (b2) node {$b_6$};
\draw[anchor = east] (b3) node {$c_3$};

\draw[anchor = east] (e3) node {$b_7 \, \,$};
\draw[anchor = west] (e4) node {$b_8$};

\draw[anchor = east] (d3) node {$b_9 \, \,$};
\draw[anchor = west] (d4) node {$\, \, b_{10}$};

\draw[anchor = west] (d5) node {$c_4$};
\draw[anchor = west] (d6) node {$c_6$};
\draw[anchor = west] (d7) node {$c_5$};

\draw[anchor = south east] (b) node {\begin{color}{red}1\end{color}};
\draw[anchor = south east] (c) node {\begin{color}{red}5\end{color}};
\draw[anchor = north west] (f1) node {\begin{color}{red}5\end{color}};
\draw[anchor = south west] (c1) node {\begin{color}{red}4\end{color}};
\draw[anchor = north] (e2) node {\begin{color}{red}3\end{color}};
\draw[anchor = south] (e1) node {\begin{color}{red}4\end{color}};
\draw[anchor = west] (b1) node {\begin{color}{red}6\end{color}};
\draw[anchor = north] (d1) node {\begin{color}{red}7\end{color}};
\draw[anchor = south] (d2) node {\begin{color}{red}2\end{color}};

\draw[anchor = south] (d) node {\begin{color}{red}17\end{color}};
\draw[anchor = west] (f4) node {\begin{color}{red}8\end{color}};
\draw[anchor = north] (f2) node {\begin{color}{red}9\end{color}};
\draw[anchor = south] (c2) node {\begin{color}{red}10\end{color}};
\draw[anchor = west] (c4) node {\begin{color}{red}11\end{color}};
\draw[anchor = west] (b2) node {\begin{color}{red}12\end{color}};
\draw[anchor = west] (e3) node {\begin{color}{red}13\end{color}};
\draw[anchor = north] (e4) node {\begin{color}{red}14\end{color}};
\draw[anchor = west] (d3) node {\begin{color}{red}15\end{color}};
\draw[anchor = north] (d4) node {\begin{color}{red}16\end{color}};

\draw[anchor = south] (d5) node {\begin{color}{red}7\end{color}};
\draw[anchor = north] (d7) node {\begin{color}{red}7\end{color}};
\draw[anchor = west] (b3) node {\begin{color}{red}6\end{color}};

\end{tikzpicture}
\end{center}
\caption{Graph $G$ with $d_G^s(v)=17$ which corresponds to a partition $A_1(v)=\{a_1^1,a_6^1,a_9^1 \}$, $A_2(v)=\{a_2^2,a_3^2,a_4^2,a_7^2\}$, $A_3(v)=\{a_5^3,a_8^3 \}$ where $X(v)=\{b_1,b_2,b_5,b_7,b_8\}$ and $Y(v)=\{b_3,b_4,b_6,b_9,b_{10}\}$.}
\label{starstopnja1}
\end{figure}
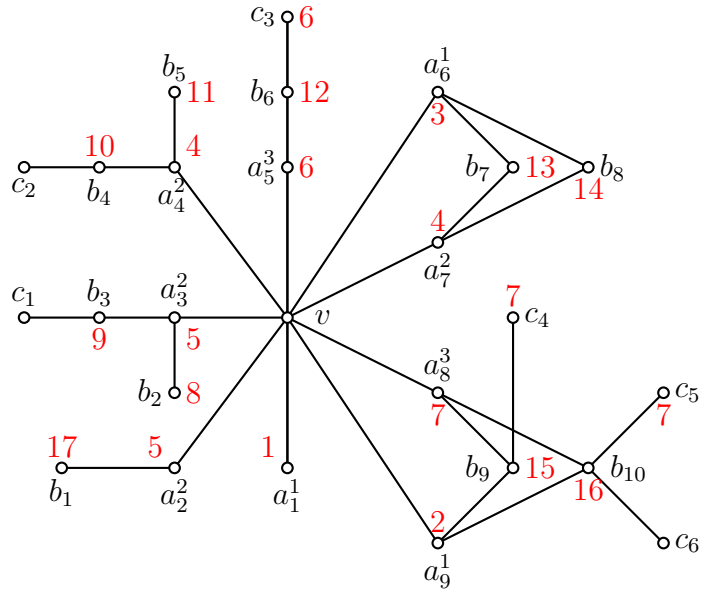

Similarly a set $A_3(v)\subseteq (N_G(v)-A_2(v))$ contains the minimum number of neighbors of $v$, such that every vertex of $Y(v)$ has a neighbor in $A_2(v)\cup A_3(v)$. As before we call $y\in Y(v)$ with $N_G(y)\cap N_G(v)=\{u\}$ a \emph{private neighbor} of $u$ in $Y(v)$. We end with $A_1(v)=N_G(v)-(A_2(v)\cup A_3(v))$. Clearly, $A_1(v),A_2(v)$ and $A_3(v)$ form a weak partition of $N_G(v)$ and are uniquely determined when $g(v)\geq 7$. A general example when $g(v)=4$ that illustrates sets $A_1(v),A_2(v),A_3(v),X(v)$ and $Y(v)$ is in Figure \ref{starstopnja1}.

\begin{theorem}\label{girth}
Let $v$ be an arbitrary vertex of a graph $G$ with $d_G(v)>1$. If $g(v)\geq 7$, then 
$$d_G^s(v)=|A_1(v)|+\left \lfloor \frac{|A_2(v)|}{2} \right \rfloor +|A_3(v)|+|X(v)|+|Y(v)|.$$
\end{theorem}

\noindent {\textbf{Proof.}}
Let $v\in V(G)$ be a vertex with $g(v)\geq 7$. Every vertex of $X(v)\cup Y(v)$ has a unique neighbor either in $A_2(v)$ or in $A_3(v)$ by the girth condition. Hence $A_1(v),A_2(v)$ and $A_3(v)$ are unique. If there exists a $0122$ blocking path $vuxy$, then $y$ is adjacent only to $u$ in $N_G(v)$ because $g(v)\geq 7$. Hence, $c(u)\neq c(y)$ which is not possible for a $0122$ blocking path. Therefore, there exists no $0122$ blocking path with respect to $v$ in $G$. Let $K=|A_1(v)|+\left \lfloor \frac{|A_2(v)|}{2} \right \rfloor +|A_3(v)|+|X(v)|+|Y(v)|$. First we show that $d_G^s(v)\geq K$. For this we describe a star coloring $c$ of $G$ where $K$ colors will be blocked for $v$. We first partition vertices of $A_2(v)$ into pairs, when $|A_2(v)|$ is an even number. If $|A_2(v)|>1$ is an odd number, then we partition $A_2(v)$ into $\frac{|A_2(v)|-3}{2}$ pairs and one triple. The case when $|A_2(v)|=1$ will be explained later. Now we color every pair (and possibly triple) with their own color $p_i$, $i \in \left [ \frac{A_2(v)}{2} \right ]$ or $i \in \left [ \frac{A_2(v)-3}{2}+1 \right ]$. In addition every vertex from $A_1(v)\cup A_3(v)\cup X(v)\cup Y(v)$ receives its own color different from $p_i$. If $|A_2(v)|=1$, then $A_1(v)\cup A_3(v)\neq\emptyset$ because $d_G(v)>1$ and we color the vertex from $A_2(v)$ with an arbitrary color that was used in $A_1(v)\cup A_3(v)$. For every vertex $y\in Y(v)$ that is adjacent to a vertex in $A_3(v)$ but not to a vertex in $A_2(v)$ we color further one neighbor of $y$, say $u_y$, from $N_3(v)$ with the same color as the neighbor $w_y$ of $y$ in $N_G(v)$: $c(u_y)=c(w_y)$. Notice that for two such vertices $y_1,y_2\in Y(v)$, we have $c(u_{y_1})=c(w_{y_1})\neq c(w_{y_2})=c(u_{y_2})$ since $w_{y_1},w_{y_2}\in A_3(v)$ and the coloring is proper even if $u_{y_1}u_{y_2}\in E(G)$. All the other vertices can be then colored arbitrarily (possibly with more colors), such that we get a star coloring of $G$. 

We claim that $K$ colors are blocked for $v$. Clearly, there are $|A_1(v)|+\left \lfloor \frac{|A_2(v)|}{2} \right \rfloor +|A_3(v)|$ different colors in $N_G(v)$ which are blocked for $v$. For a vertex $y\in Y(v)$ with a neighbor in $A_3(v)$, the color $c(y)$ is blocked by the $0123$ blocking path $vw_yyu_y$. Let now $x\in X(v)$ and let $x$ be a private neighbor of $w_x\in A_2$. By the construction of the coloring, there exists a vertex $z\in N_G(v)$ different than $w_x$ with $c(z)=c(w_x)$. The path $zvw_xx$ is $1012$ blocking for $v$ so $c(x)$ is also blocked for $v$. All in all we have $K$ blocked colors for $v$ and  $d_G^s(v)\geq K$ follows.

Assume now, on the way to a contradiction, that $d_G^s(v)>K$. The color of every neighbor of $v$ is blocked for $v$. Since $g(v)\geq 7$, every $1012$ and every $0123$ blocking path can provide exactly one additional color that is blocked for $v$. In other words, we have at most $|X(v)|+|Y(v)|$ additional colors that are not in $N_G(v)$ that can be blocked for $v$. Let $D$ denote the set of all vertices in $A_2(v)$ whose color is different from the colors of all other neighbors of $v$. By the pigeon hole principle $D\neq\emptyset$. Every $z\in D$ has at least one private neighbor $x_z\in X(v)$ and $c(x_z)$ is not blocked for $v$. So we have at most 
$$|A_1(v)|+\left \lfloor \frac{|A_2(v)-D|}{2} \right \rfloor +|D|+|A_3(v)|+|X(v)|-|D|+|Y(v)|$$
colors that are blocked for $v$, which is $\leq K$, a contradiction. 
\qed\bigskip

The condition in above theorem that $d_G(v)>1$ is necessary which one can observe from a star $K_{1,t}$, $t\geq 3$. For a leaf $v\in V(K_{1,t})$ it holds that $d_{K_{1,t}}^s(v)=1$. However, the theorem states that $d_{K_{1,t}}^s(v)=t-1$, which is greater than $1$ for every $t\geq 3$.

If $g(v)<7$, then we can have less blocked colors in $X(v)$ or in $Y(v)$. The first one can occur when $x\in X(v)$ is a private neighbor of $u\in A_2(v)$ and $u$ cannot be colored with the same color as some other vertex from $N_G(v)$. In the next result we omit this possibility, but we need to take care of the second option. As explained earlier, two $0123$ blocking paths with the same last vertex can force less blocked colors in $Y(v)$. Therefore we introduce $Y_1(v)$ as a subset of $Y(v)$ of maximum cardinality such that there exists a perfect matching between $Y_1(v)$ and some subset of $N_3(v)$.

\begin{theorem}\label{matching}
Let $v$ be an arbitrary vertex of a graph $G$ that does not belong to a five cycle. If there exist $A_2(v)$ with a perfect matching $M$ in $\overline{G}[A_2(v)]$, then 
$$d_G^s(v)=|A_1(v)|+\frac{|A_2(v)|}{2}+|A_3(v)|+|X(v)|+|Y_1(v)|.$$
\end{theorem}

\noindent {\textbf{Proof.}}
Let $v\in V(G)$ be a vertex and let $A_2(v)$ contain a perfect matching $M$ in $\overline{G}[A_2(v)]$. As in the proof of Theorem \ref{girth}, no $0122$ blocking path exists here as well by the same reason as in the mentioned proof or since $v$ does not belong to a five cycle. Let $K=|A_1(v)|+\frac{|A_2(v)|}{2}+|A_3(v)|+|X(v)|+|Y_1(v)|$. First we show that $d_G^s(v)\geq K$ by a star coloring $c$ of $G$ where $K$ colors are blocked for $v$. For this we first partition vertices of $A_2(v)$ into pairs according to the pairs of $M$. Now we color every pair with their own color. In addition, every vertex in $A_1(v)\cup A_3(v)\cup X(v)\cup Y_1(v)$ is assigned a unique colors different from the colors of all the pairs. Let $M'$ be a perfect matching between $Y_1(v)$ and a subset of $N_3(v)$. So, for every $y\in Y_1(v)$ there exists $u_y\in N_3(v)$ according to $M'$ and we color $u_y$ with the same color as the neighbor $w_y\in N_G(v)$ of $y$: $c(u_y)=c(w_y)$. All the other vertices can then be colored arbitrarily (with possibly more colors), such that we get a star coloring of $G$. 

We claim that $K$ colors are blocked for $v$. Clearly, there are $|A_1(v)|+\frac{|A_2(v)|}{2}+|A_3(v)|$ different colors in $N_G(v)$ which are blocked for $v$. For a vertex $y\in Y_1(v)$ the color $c(y)$ is blocked by the $0123$ blocking path $vw_yyu_y$. Let now $x\in X(v)$ and let $x$ be a private neighbor of $w_x\in A_2(v)$.  There exists $z\in A_2(v)$, such that $(w_x,z)\in M$, meaning that $c(z)=c(w_x)$. The path $zvw_xx$ is $1012$ blocking for $v$ so $c(x)$ is also blocked for $v$. All in all we have $K$ blocked colors for $v$ and $d_G^s(v)\geq K$ follows.

Assume now, on the way to a contradiction, that $d_G^s(v)>K$. The color of every neighbor of $v$ is blocked for $v$. As mentioned, no $0122$ blocking path exists for $v$. Every $0123$ or $1012$ blocking path blocks at most one additional color. Moreover, every $0123$ blocking path $vwyu$ that blocks one additional color beside $c(w)$, has $c(w)=c(u)$.

Altogether we have at most $|X(v)|+|Y_1(v)|$ additional colors that are not in $N_G(V)$ that can be blocked for $v$. Let $D$ denote the set of all vertices in $A_2(v)$ whose color is different from the colors of all other neighbors of $v$. By the pigeon hole principle $D\neq\emptyset$. Every $z\in D$ has at least one private neighbor $x_z\in X$ and $c(x_z)$ is not blocked for $v$. So we have at most 
$$|A_1(v)|+\frac{|A_2(v)-D|}{2}+|D|+|A_3(v)|+|X(v)|-|D|+|Y_1(v)|$$
colors that are blocked for $v$, which is $\leq K$, a contradiction. 
\qed\bigskip

\begin{figure}[ht!]
\begin{center}
\begin{tikzpicture}[xscale=.7, yscale=.7,style=thick,x=0.8cm,y=0.8cm]
\def\vr{2pt} 

\path (0,0) coordinate (a);
\path (2,0) coordinate (b);
\path (4,0) coordinate (c);
\path (6,0) coordinate (v);

\path (8,2) coordinate (b1);
\path (8,-2) coordinate (b2);
\path (10,2) coordinate (b3);
\path (10,-2) coordinate (b4);

\path (5,2) coordinate (c1);
\path (5,-2) coordinate (c2);
\path (2.5,3) coordinate (c3);
\path (2.5,-3) coordinate (c4);

\path (6.5,2) coordinate (d1);
\path (6.5,-2) coordinate (d2);
\path (4.5,3) coordinate (d3);
\path (4.5,-3) coordinate (d4);


\draw (a) -- (b) -- (c) -- (v) -- (b1) -- (b3);
\draw (c4) -- (c2) -- (v);
\draw (b4) -- (b2) -- (v) -- (d2) -- (d4) -- (c4) -- (a) -- (c3) -- (d3) -- (d1) -- (v) -- (c1) -- (c3);

\draw (b) [fill=white] circle (\vr);
\draw (c) [fill=white] circle (\vr);
\draw (a) [fill=white] circle (\vr);
\draw (v) [fill=white] circle (\vr);
\draw (b1) [fill=white] circle (\vr);
\draw (b2) [fill=white] circle (\vr);
\draw (b3) [fill=white] circle (\vr);
\draw (b4) [fill=white] circle (\vr);

\draw (c1) [fill=white] circle (\vr);
\draw (c2) [fill=white] circle (\vr);
\draw (c3) [fill=white] circle (\vr);
\draw (c4) [fill=white] circle (\vr);

\draw (d1) [fill=white] circle (\vr);
\draw (d2) [fill=white] circle (\vr);
\draw (d3) [fill=white] circle (\vr);
\draw (d4) [fill=white] circle (\vr);

\draw[anchor = west] (v) node {$v$};
\draw[anchor = south] (b) node {$5$};
\draw[anchor = south] (c) node {$4$};
\draw[anchor = south east] (a) node[xshift=1mm]  {$4$};

\draw[anchor = north east] (c1) node[xshift=1mm] {$7$};
\draw[anchor = south east] (c2) node[xshift=1mm] {$6$};
\draw[anchor = south] (c3) node {$11$};
\draw[anchor = north] (c4) node {$8$};

\draw[anchor = south] (b1) node {$1$};
\draw[anchor = north] (b2) node {$1$};
\draw[anchor = south] (b3) node {$3$};
\draw[anchor = north] (b4) node {$2$};

\draw[anchor = south] (d1) node {$11$};
\draw[anchor = north] (d2) node {$8$};
\draw[anchor = south] (d3) node {$10$};
\draw[anchor = north] (d4) node {$9$};

\end{tikzpicture}
\end{center}
\caption{Vertex $v$ with $d_G^s(v)=11$ where colors $9$ and $10$ are blocked by $0122$ blocking paths.}
\label{starstopnja}
\end{figure}
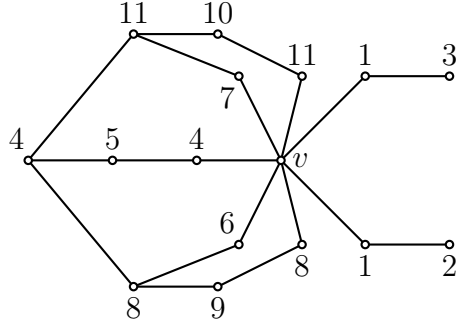

We end this section with an example in Figure \ref{starstopnja} where $0122$ blocking paths bring some additional blocked colors for $v$. 

\section{Two upper bounds of $S_b(G)$}\label{secupperbound} 

In this section we define a natural upper bound for $S_b(G)$. The main idea is that we need a sufficiently large number of vertices with high star degree as candidates for star b-vertices which is justified by Corollary \ref{cor}. On the other hand, it is a direct generalization of $m$-degree that is a natural upper bound for $\varphi(G)$, see \cite{IrMa}, and of $m_a$-degree that is an upper bound of acyclic b-chromatic number, see \cite{AnCP}.   
Let the vertices $v_1,\dots,v_{n_G}$ of $G$ be ordered by their non-increasing star degree $d_G^s(v_1)\geq \cdots\geq d_G^s(v_{n_G})$. We define an $m_s$-degree of a graph $G$, denoted by $m_s(G)$, as 
$$m_s(G)=\max \{i: i-1\leq d_G^s(v_i)\}.$$

\begin{theorem}\label{madegree}
For any graph $G$ we have $S_b(G)\leq m_s(G)$.
\end{theorem}

\noindent {\textbf{Proof.}}
Suppose on the contrary that $S_b(G)>m_s(G)$. Let $c:V(G)\rightarrow [S_b(G)]$ be an $S_b(G)$-coloring of $G$. From the definition of $m_s(G)$ it follows that there exist $m_s(G)$ vertices having star degree at least $m_s(G)-1$. In other words, this means that $G$ contains $|V(G)|-m_s(G)$ vertices whose star degree is at most $m_s(G)-1$. Since we assumed that $S_b(G)>m_s(G)$, there exists at least one color $\ell\in [k]$ such that, for every vertex $v$ of color $\ell$ it holds that $d_G^s(v)\leq m_s(G)-1$. Therefore, for every vertex $v\in V(G)$ of color $\ell$ we have $d_G^s(v)< S_b(G)-1$, meaning that no vertex of color $\ell$ is a star b-vertex of $c$, a contradiction by Corollary \ref{cor}.
\qed\bigskip

We can use $m_s$-degree to obtain an upper bound for $S_b(G)$ with respect to $\Delta(G)$ as follows.

\begin{proposition}\label{zgornjameja}
For any graph $G$ we have $m_s(G)\leq \Delta(G)^2+1$.
\end{proposition}

\noindent {\textbf{Proof.}} 
For any vertex $v$ at most $\Delta(G)$ colors can be blocked in $N_G(v)$. In addition, every blocking path blocks at most one new color. So, for every neighbor of $v$, this gives at most $\Delta(G)-1$ new colors. Hence,   
$$d_G^s(v)\leq \Delta(G)+\Delta(G)(\Delta (G) -1)=\Delta(G)^2.$$
Consequently, $m_s(G)\leq d_G^S(v)+1\leq \Delta(G)^2+1$ by Theorem \ref{madegree}.
\qed\bigskip

\begin{corollary}\label{zgornjameja1}
For any graph $G$ we have $S_b(G)\leq \Delta(G)^2+1$.
\end{corollary}

The bounds from the last three results are sharp even for some trees as shown next.

\begin{theorem}\label{DeltaSquaredExtremal}
There exists a family of trees $T_1, T_2, \dots$ such that $S_b(T_n)=\Delta(T_n)^2+1$. 
\end{theorem}

\noindent {\textbf{Proof.}}
For any positive integer $n$, we are going to define $n^2+1$ copies $H_{n,i}$, $i\in[n^2+1]$, of the same tree $H_n$.  Furthermore, we define the tree $T_n$ using the mentioned copies, making it a member of the desired family. The vertices and edges of $H_{n,i}$ are defined as follows: 
\begin{align*}
V(H_{n,i})&=\{v^i\}\cup\{x^i_j:j\in[n]\}\cup\{y^i_k:k\in[n(n-1)]\}\cup\{z^i_{\ell}:\ell\in[n(n-1)]\},\\
E(H_{n,i})&=\{v^ix^i_j:j\in[n]\}\\
&\cup\{x^i_jy^i_k: j\in [n], k\in[(n-1)(j-1)+1,(n-1)j]\}\\ 
&\cup\{y^i_kz^i_k: k\in [n(n-1)]\}.
\end{align*}

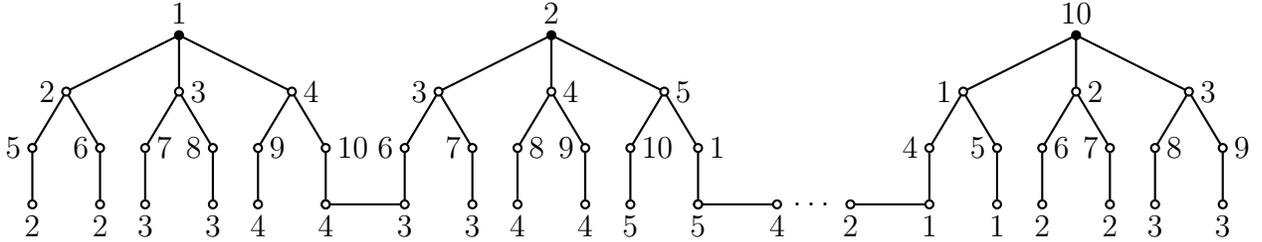
\begin{figure}[ht!]
\begin{center}
\begin{tikzpicture}[xscale=0.75, yscale=0.75,style=thick,x=1cm,y=1cm]
\def\vr{2pt} 

\path (2,3) coordinate (v);
\path (0,2) coordinate (a);
\path (2,2) coordinate (b);
\path (4,2) coordinate (c);

\path (-0.6,1) coordinate (a1);
\path (0.6,1) coordinate (a2);
\path (1.4,1) coordinate (b1);
\path (2.6,1) coordinate (b2);
\path (3.4,1) coordinate (c1);
\path (4.6,1) coordinate (c2);

\path (-0.6,0) coordinate (a3);
\path (0.6,0) coordinate (a4);
\path (1.4,0) coordinate (b3);
\path (2.6,0) coordinate (b4);
\path (3.4,0) coordinate (c3);
\path (4.6,0) coordinate (c4);


\draw (v) -- (a) -- (a1) -- (a3);
\draw (a) -- (a2) -- (a4);
\draw (v) -- (b) -- (b1) -- (b3);
\draw (b) -- (b2) -- (b4);
\draw (v) -- (c) -- (c1) -- (c3);
\draw (c) -- (c2) -- (c4);

\draw (v) [fill=black] circle (\vr);
\draw (a) [fill=white] circle (\vr);
\draw (b) [fill=white] circle (\vr);
\draw (c) [fill=white] circle (\vr);
\draw (a1) [fill=white] circle (\vr);
\draw (a2) [fill=white] circle (\vr);
\draw (a3) [fill=white] circle (\vr);
\draw (a4) [fill=white] circle (\vr);
\draw (b1) [fill=white] circle (\vr);
\draw (b2) [fill=white] circle (\vr);
\draw (b3) [fill=white] circle (\vr);
\draw (b4) [fill=white] circle (\vr);
\draw (c1) [fill=white] circle (\vr);
\draw (c2) [fill=white] circle (\vr);
\draw (c3) [fill=white] circle (\vr);
\draw (c4) [fill=white] circle (\vr);

\draw[anchor = south] (v) node {1};
\draw[anchor = east] (a) node {2};
\draw[anchor = west] (b) node {3};
\draw[anchor = west] (c) node {4};

\draw[anchor = east] (a1) node {5};
\draw[anchor = east] (a2) node {6};
\draw[anchor = west] (b1) node {7};
\draw[anchor = east] (b2) node {8};
\draw[anchor = west] (c1) node {9};
\draw[anchor = west] (c2) node {10};

\draw[anchor = north] (a3) node {2};
\draw[anchor = north] (a4) node {2};
\draw[anchor = north] (b3) node {3};
\draw[anchor = north] (b4) node {3};
\draw[anchor = north] (c3) node {4};
\draw[anchor = north] (c4) node {4};

\path (8.6,3) coordinate (u);
\path (6.6,2) coordinate (d);
\path (8.6,2) coordinate (e);
\path (10.6,2) coordinate (f);

\path (6,1) coordinate (d1);
\path (7.2,1) coordinate (d2);
\path (8,1) coordinate (e1);
\path (9.2,1) coordinate (e2);
\path (10,1) coordinate (f1);
\path (11.2,1) coordinate (f2);

\path (6,0) coordinate (d3);
\path (7.2,0) coordinate (d4);
\path (8,0) coordinate (e3);
\path (9.2,0) coordinate (e4);
\path (10,0) coordinate (f3);
\path (11.2,0) coordinate (f4);

\draw (c4) -- (d3);
\draw (c4) [fill=white] circle (\vr);

\draw (u) -- (d) -- (d1) -- (d3);
\draw (d) -- (d2) -- (d4);
\draw (u) -- (e) -- (e1) -- (e3);
\draw (e) -- (e2) -- (e4);
\draw (u) -- (f) -- (f1) -- (f3);
\draw (f) -- (f2) -- (f4);

\draw (u) [fill=black] circle (\vr);
\draw (d) [fill=white] circle (\vr);
\draw (e) [fill=white] circle (\vr);
\draw (f) [fill=white] circle (\vr);
\draw (d1) [fill=white] circle (\vr);
\draw (d2) [fill=white] circle (\vr);
\draw (d3) [fill=white] circle (\vr);
\draw (d4) [fill=white] circle (\vr);
\draw (e1) [fill=white] circle (\vr);
\draw (e2) [fill=white] circle (\vr);
\draw (e3) [fill=white] circle (\vr);
\draw (e4) [fill=white] circle (\vr);
\draw (f1) [fill=white] circle (\vr);
\draw (f2) [fill=white] circle (\vr);
\draw (f3) [fill=white] circle (\vr);
\draw (f4) [fill=white] circle (\vr);

\draw[anchor = south] (u) node {2};
\draw[anchor = east] (d) node {3};
\draw[anchor = west] (e) node {4};
\draw[anchor = west] (f) node {5};

\draw[anchor = east] (d1) node {6};
\draw[anchor = east] (d2) node {7};
\draw[anchor = west] (e1) node {8};
\draw[anchor = east] (e2) node {9};
\draw[anchor = west] (f1) node {10};
\draw[anchor = west] (f2) node {1};

\draw[anchor = north] (d3) node {3};
\draw[anchor = north] (d4) node {3};
\draw[anchor = north] (e3) node {4};
\draw[anchor = north] (e4) node {4};
\draw[anchor = north] (f3) node {5};
\draw[anchor = north] (f4) node {5};

\path (12.6,0) coordinate (z1);
\path (13.9,0) coordinate (z2);
\node at (13.25,0) (n) {$\dots$};
\draw (f4) -- (z1);
\path (15.3,0) coordinate (g3);
\draw (z2) -- (g3);
\draw (f4) [fill=white] circle (\vr);
\draw (z1) [fill=white] circle (\vr);
\draw (z2) [fill=white] circle (\vr);

\draw[anchor = north] (z1) node {4};
\draw[anchor = north] (z2) node {2};

\path (17.9,3) coordinate (w);
\path (15.9,2) coordinate (g);
\path (17.9,2) coordinate (h);
\path (19.9,2) coordinate (i);

\path (15.3,1) coordinate (g1);
\path (16.5,1) coordinate (g2);
\path (17.3,1) coordinate (h1);
\path (18.5,1) coordinate (h2);
\path (19.3,1) coordinate (i1);
\path (20.5,1) coordinate (i2);

\path (16.5,0) coordinate (g4);
\path (17.3,0) coordinate (h3);
\path (18.5,0) coordinate (h4);
\path (19.3,0) coordinate (i3);
\path (20.5,0) coordinate (i4);

\draw (w) -- (g) -- (g1) -- (g3);
\draw (g) -- (g2) -- (g4);
\draw (w) -- (h) -- (h1) -- (h3);
\draw (h) -- (h2) -- (h4);
\draw (w) -- (i) -- (i1) -- (i3);
\draw (i) -- (i2) -- (i4);

\draw (w) [fill=black] circle (\vr);
\draw (g) [fill=white] circle (\vr);
\draw (h) [fill=white] circle (\vr);
\draw (i) [fill=white] circle (\vr);
\draw (g1) [fill=white] circle (\vr);
\draw (g2) [fill=white] circle (\vr);
\draw (g3) [fill=white] circle (\vr);
\draw (g4) [fill=white] circle (\vr);
\draw (h1) [fill=white] circle (\vr);
\draw (h2) [fill=white] circle (\vr);
\draw (h3) [fill=white] circle (\vr);
\draw (h4) [fill=white] circle (\vr);
\draw (i1) [fill=white] circle (\vr);
\draw (i2) [fill=white] circle (\vr);
\draw (i3) [fill=white] circle (\vr);
\draw (i4) [fill=white] circle (\vr);

\draw[anchor = south] (w) node {10};
\draw[anchor = east] (g) node {1};
\draw[anchor = west] (h) node {2};
\draw[anchor = west] (i) node {3};

\draw[anchor = east] (g1) node {4};
\draw[anchor = east] (g2) node {5};
\draw[anchor = west] (h1) node {6};
\draw[anchor = east] (h2) node {7};
\draw[anchor = west] (i1) node {8};
\draw[anchor = west] (i2) node {9};

\draw[anchor = north] (g3) node {1};
\draw[anchor = north] (g4) node {1};
\draw[anchor = north] (h3) node {2};
\draw[anchor = north] (h4) node {2};
\draw[anchor = north] (i3) node {3};
\draw[anchor = north] (i4) node {3};

\end{tikzpicture}
\end{center}
\caption{Graph $T_3$ and its star b-coloring with $10$ colors.}
\label{zgornjamejaslika}
\end{figure}

First, let $T_1=H_{1,1}=P_2$. Then clearly $S_b(T_1)=2=(\Delta(T_1))^2+1$. Additionally, for any $v\in V(T_1)$ it holds $d_{T_1}^s(v)=1$, thus $m_s(T_1)=2$. Now, let $n\geq 2$. The graph $T_n$ is obtained from $n^2+1$ graphs $H_{n,i}$, $i\in[n^2+1]$, by adding the set of edges $\{z^i_{n(n-1)}z^{i+1}_1:i\in[n^2]\}$. Notice that $T_2\cong P_{35}$. Graph $T_3$ and its 10-star b-coloring are presented in Figure \ref{zgornjamejaslika}, where the star b-vertices of each color class are colored in black.
Note that $T_n$ has $n^2+1$ vertices $v^i$, $i\in[n^2+1]$, of degree $n$, $(n^2+1)n$ vertices $x_j^i$, $i\in[n^2+1]$, $j\in[n]$ of degree $n$, $(n^2+1)n(n-1)$ vertices $y^i_k$, $i\in[n^2+1]$, $k\in[n(n-1)]$, of degree $2$ and $(n^2+1)n(n-1)$ vertices $z^i_{\ell}$, $i\in[n^2+1]$, $\ell\in[n(n-1)]$, where $2n^2$ of them have degree $2$ and all the others have degree $1$. 
It follows that $\Delta(T_n)=n$ and $S_b(G_n)\leq \Delta(T_n)^2+1=n^2+1$ by Corollary \ref{zgornjameja1}. 

To complete the proof, we construct a coloring $c$ that uses elements of the additive group $\mathbb{Z}_{n^2+1}$ as colors:
\begin{align*}
    c(v^i)&=i, &&i\in[n^2+1],\\
    c(x^i_j)&= i+j, &&j\in[n],\\
    c(y^i_k)&=i+n+k, &&k\in[n(n-1)],\\
    c(z^i_{\ell})&=c(x^i_j), &&j\in[n], \ell\in[(n-1)(j-1)+1,(n-1)j].\\
\end{align*}
It is easy to see that $c$ is a star coloring of $T_n$ in which every $v^i$, $i\in[n^2+1]$, is a star b-vertex of color $i$. Namely, exactly $n^2$ colors are blocked for every vertex $v^i$ ($n$ colors by its neighbors $x_j^i$ and $n(n-1)$ remaining colors by paths $v^ix_j^iy_k^iz_k^i$, by which any of the potential $n(n-1)$ recolorings of vertex $v^i$ is unachievable). Because a vertex $v^i$ exists in every color class of $c$, we infer that $c$ is a star b-coloring of $T_n$ by Corollary \ref{cor}. Consequently, we conclude that $S_b(T_n)\geq n^2+1=(\Delta(T_n))^2+1$ and the equality follows. \qed\bigskip

One could also apply Theorem \ref{girth} to obtain $d^s(v^i)$, $i\in[n^2+1]$, of trees $T_n$ because $g(v^i)=\infty$. For this, the weak partition $A_0=A_1=\emptyset$, $A_2=\{x^i_j:j\in[n]\}$ of $N_{T_n}(v^i)$, we have $|A_0|=0$, $\left \lfloor \frac{|A_1|}{2} \right \rfloor=0$, $|X(v^i)|=0$, $|A_2|=n$ and $|Y(v^i)|=n(n-1)$. Consequently, $d_{T_n}^s(v^i)\geq 0+0+0+n+n(n-1)=n^2=(\Delta(T_n))^2$.


\section{Some exact results and relations between $S_b(G)$ and other parameters}

To demonstrate that the star b-chromatic number of a tree can be intriguing, we next consider paths. 

\begin{theorem}\label{poti}
Let $P_n$ be a path on $n$ vertices. Then 
$$S_b(P_n)=
		\begin{cases}
			1 &; \, n=1 \\
			2 &; \, 2 \leq n \leq 3 \\
			3 &; \, 4 \leq n \leq 7 \\
			4 &; \, 8 \leq n \leq 22 \\
			5 &; \, n \geq 23. \\
		\end{cases}$$
\end{theorem}

\noindent {\textbf{Proof.}}
Let $P_n=v_1\dots v_n$ be a path on $n$ vertices. Obviously, $S_b(P_1)=1$ and $S_b(P_2)=2$. Also at least 2 colors are needed to obtain an $S_b(P_3)$-coloring. Suppose there exists an $S_b(P_3)$-coloring with 3 colors. That means that every vertex receives its own color. But then any of the leaves can be recolored, which is a contradiction. Thus, $S_b(P_3)=2$. 

To prove that $S_b(P_4)\geq 3$ we construct a 3-star b-coloring $1-2-1-3$ of $P_4$. Further, $S_b(P_k)\geq 3$ for $k\in [5,7]$ by extending the mention coloring of $P_4$ by $2$ for $P_5$, by $2,1$ for $P_6$ and by $2,1,3$ for $P_7$. 

Next, we show that there does not exist an $S_b(P_7)$-coloring with 4 colors. Suppose on the contrary that $c$ is an $S_b(P_7)$-coloring using the colors $[4]$. This means that every color class of $c$ has at least one star b-vertex. Clearly, $d^s(v_1)=d^s(v_7)=2$ and $v_1$ and $v_7$ are not star b-vertices for $c$. At least one of $v_2$ and $v_6$, say $v_2$, must be a star b-vertex for color $j$ since we need four star b-vertices. Let first $c(v_{1})=c(v_{3})=\ell$, $\ell\in[4]-\{j \}$. Only one color from $[4]-\{j,\ell\}$ can be blocked for $v_2$ in this configuration, so this is not possible. Thus, we may assume that $k=c(v_{1})\neq c(v_{3})=\ell$ where $[4]=\{j,k,\ell,m\}$. Color $m$ can now be blocked for $v_2$ only by setting $c(v_{4})=m$ and $c(v_{5})=\ell$. Now, to obtain a star b-vertex of color $k$ we are forced to set $c(v_{6})=k$ and $c(v_7)=j$. However, color $\ell$ is without a star b-vertex, a contradiction. Thus, an $S_b(P_7)$-coloring using four colors does not exists and $S_b(P_7)\leq 3$. By similar reasons as for $P_7$ there is also no $S_b(P_k)$-coloring with four colors for any $k\in [4,6]$. We conclude that $S_b(P_n)=3$ for every $n\in [4,7]$.

To prove that $S_b(P_{8})\geq 4$ we construct a star b-coloring $4-3-2-1-2-4-3-4$ of $P_8$. For $k\in [9,22]$ we extend the mentioned coloring by repeating colors $1,2,3$ as long as needed to obtain a star b-coloring with four colors of $P_k$ for every $k\in [9,22]$. So, $S_b(P_k)\geq 4$ for every $k\in [8,22]$.

Further, we assume there exists an $S_b(P_{n})$-coloring, $n\in [8,22]$, with 5 colors where $[5]=\{j,k,\ell,m,p\}$. The only way for vertex $v_i$ to be a star b-vertex of color $j$ is to block two colors, say $k$ and $\ell$, as neighbors and the other two colors $m$ and $p$, by $0123$ blocking paths. This means that $i\in [4,19]$ and without loss of generality $c(v_{i-1})=k=c(v_{i-3})$, $c(v_{i+1})=\ell=c(v_{i+3})$, $c(v_{i-2})=m$ and $c(v_{i+2})=p$. Moreover, $v_a$, $a\in [i-3,i+3]-\{i\}$, is not a star b-vertex. Consequently, $2\cdot 3+1=7$ vertices are needed to obtain one star b-vertex (with the middle vertex becoming a star b-vertex). Here $v_{i-3},v_{i-2},v_{i-1}$ can influence $v_{i-4}$ to become a star b-vertex and symmetrically $v_{i+1},v_{i+2},v_{i+3}$ can influence $v_{i+4}$. Altogether, we need at least $6\cdot 3+5=23$ vertices to obtain a star b-coloring with 5 colors, a contradiction. Thus, $S_b(P_{n})=4$ for any $n\in [8,22]$.

Next we consider $P_{23}$. It is straightforward to check that 
\begin{equation}\label{coloring}
1-4-1-5-2-3-2-4-1-5-1-2-4-3-4-1-2-5-2-3-1-4-1
\end{equation}
is a star b-coloring of $P_{23}$ with 5 colors. Thus, $S_b(P_{23})\geq 5$. For $S_b(P_n)$, $n\geq 24$, we just extend the mentioned coloring with repeating the pattern $2,3,4$ as long as needed. Hence, $S_b(P_{n})\geq 5$ for any $n\geq 23$. By Corollary \ref{zgornjameja1} we have $S_b(P_{n})\leq 5$ for any $n\geq 23$ and equality follows. \qed \bigskip

The values of $S_b(G)$ for cycles are similar to those for paths, with the only difference being in the values of $n$.

\begin{theorem}\label{cikli}
Let $C_n$ be a cycle on $n\geq 3$ vertices. Then 
$$S_b(C_n)=
		\begin{cases}
			3 &; \, n \leq 8 \wedge n\neq 5\\
			4 &; \, n=5 \vee (9 \leq n \leq 23 \wedge n\neq 20) \\
			5 &; \, n=20 \vee n\geq 24. \\
		\end{cases}$$
\end{theorem}

\noindent {\textbf{Proof.}}
For $C_{20}$ notice that the coloring (\ref{coloring}) has the same first and last three colors and if we amalgamate them, we obtain a star b-coloring of $C_{20}$ with five colors. For $C_{24}$ observe the coloring (\ref{coloring}) again where color 5 is added to the non-colored vertex. For longer cycles just extend the mentioned coloring by repeating the pattern $2,3,5$ as long as needed. By Corollary \ref{zgornjameja1} we have $S_b(C_{n})\leq 5$ for any cycle and equality of the last line follows. 

For $C_n$, $n<20$ we can use a similar explanation as in paragraph above (\ref{coloring}) to see that there are not enough vertices to ensure a star b-coloring with five colors. Hence, $S_b(C_n)\leq 4$. For $n\in [21,23]$ we further analyze the coloring (\ref{coloring}). After the first star b-vertex is settled in (\ref{coloring}), the next star b-vertex is assigned one of the two starting colors, and the third star b-vertex is assigned one of the fourth and fifth color. Once this is chosen, everything else is forced. In particular last three vertices are colored with the same two colors as first three. So, if we have a cycle $C_{23}$ and try to color it with five star b-vertices, then we get a non-proper coloring or a bi-colored $P_4$, which is not a star coloring. Thus $S_b(C_{23})\leq 4$ and by the same reason also $S_b(C_n)\leq 4$ for $n\in \{21,22\}$. 

To prove that $S_b(C_{9})\geq 4$ observe a 4-star b-coloring $1-4-1-2-3-2-4-1-3$. For longer cycles $C_n$ where $10\leq n\leq 19$ or $21\leq n\leq 23$ just extend the mentioned coloring by repeating either the pattern $4,3,2$ when $n\equiv 0 ({\rm mod}\ 3)$ or the pattern $2,3,4$ otherwise as long as needed.    

For the rest, we have $S_b(C_3)=S_b(K_3)=3$ by Proposition \ref{order}. By the same proposition we have $S_b(C_4)<4$ and $S_b(C_4)=3$ follows by a star b-coloring $1-2-1-3$. Similarly we have $S_b(C_5)<5$ by Proposition \ref{order} and $S_b(C_5)=4$ follows by a star b-coloring $1-2-1-3-4$ (there is a $P_4$-system for color $1$). 

Let $C_n=v_1\dots v_nv_1$.  
To prove that $S_b(C_n)\leq 3$, $n\in [6,8]$, suppose on the contrary that $c$ is an $S_b(C_n)$-coloring using the set of colors $[4]$. Then, there exists at least one star b-vertex of every color class of $c$. Without loss of generality suppose that $v_i \in V(C_n)$ is a star b-vertex of the color $j\in[4]$. The neighbors of $v_i$ can have the same color or two distinct color different from $j$. In the beginning we consider the first possibility where $c(v_{i-1})=c(v_{i+1})=\ell$, $\ell\in[4]\setminus \{j \}$. The remaining two colors $m, k\in [4]\setminus \{j,\ell\}$ can now be blocked for $v_i$ only by blocking path $1012$ meaning that $c(v_{i-2})=m$ and $c(v_{i+2})=k$. But now the obtained coloring cannot be extended in such a way that there would exist the star b-vertices for all colors $m$, $k$ and $\ell$. Next consider the second possibility, where the neighbors of $v_i$ receive two distinct colors, say $c(v_{i-1})=k$ and $c(v_{i+1})=\ell$. To block the remaining color $m\in[4] \setminus \{j,k,\ell \}$ for $v_i$ we have to have blocking path $0123$. Meaning that $c(v_{i-2})=m$ and $c(v_{i-3})=k$ or $c(v_{i+2})=m$ and $c(v_{i+3})=\ell$, say the later one. For $v_{i+2}$ to be a star b-vertex of color $m$ we have to block the remaining color $k$ by blocking path $1012$. This is achieved by setting $c(v_{i+4})=k$. Again, the obtained coloring cannot be extended in such a way that there would exist star b-vertices of both colors $k$ and $\ell$, which means that there does not exist a 4-star b-coloring of $C_n$. We conclude that $S_b(C_n)=3$ for any $n\in[6,8]$ by star b-colorings $1-2-3-1-2-3$ for $C_6$ and adding $2$ for $C_7$ and coloring $1-2-1-3-1-2-1-3$ for $C_8$. \qed \bigskip

Next, we show that $m_s$-degree can behave unfavorably with respect to $S_b(G)$. Again, trees suffice to demonstrate that.

\begin{theorem}\label{S-m-Infinite}
There exists an infinite family of trees $T_1, T_2, \dots$ such that $(m_s(T_n)-S_b(T_n))\rightarrow\infty$ as $n\rightarrow\infty$. 
\end{theorem}

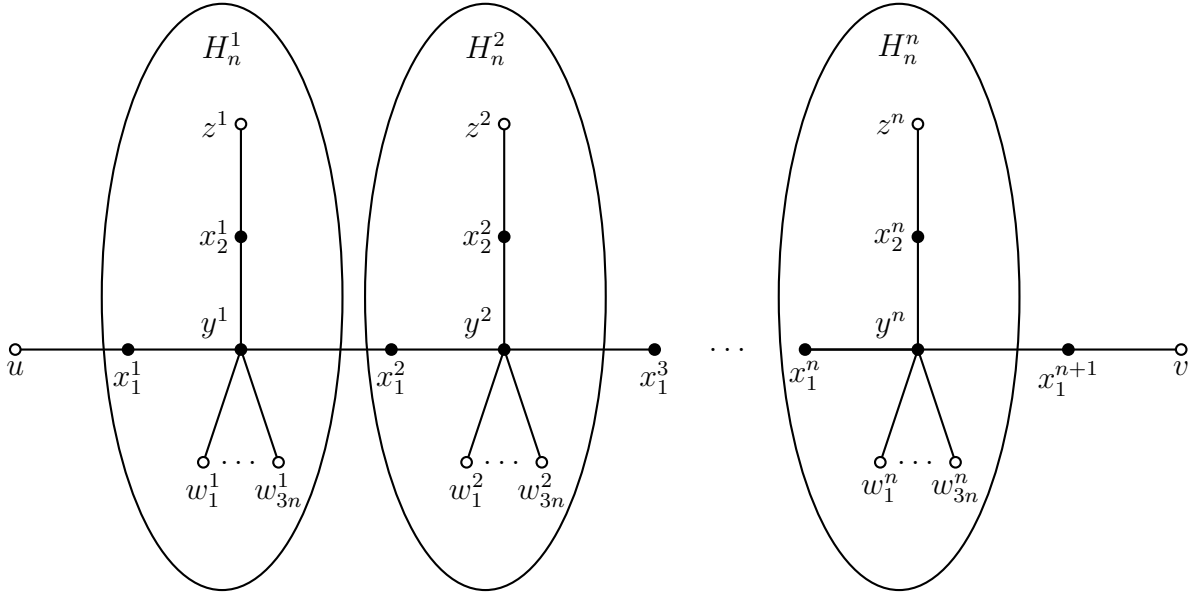
\begin{figure}[ht!]
\begin{center}
\begin{tikzpicture}[xscale=1, yscale=1,style=thick,x=1cm,y=1cm]
\def\vr{2pt} 

\path (0,0) coordinate (u);

\path (1.5,0) coordinate (x11);
\path (3,0) coordinate (y1);
\path (3,1.5) coordinate (x21);
\path (3,3) coordinate (z1);
\path (2.5,-1.5) coordinate (w11);
\path (3.5,-1.5) coordinate (w3n1);

\draw (u) -- (x11);
\draw (x11) -- (y1) -- (x21) -- (z1);
\draw (y1) -- (w11);
\draw (y1) -- (w3n1);

\draw (u) [fill=white] circle (\vr);

\draw (x11) [fill=black] circle (\vr);
\draw (y1) [fill=black] circle (\vr);
\draw (x21) [fill=black] circle (\vr);
\draw (z1) [fill=white] circle (\vr);
\draw (w11) [fill=white] circle (\vr);
\node at (3,-1.5) (n) {$\dots$};
\draw (w3n1) [fill=white] circle (\vr);

\draw[anchor = north] (u) node {$u$};

\draw[anchor = north] (x11) node {$x_1^1$};
\draw[anchor = south east] (y1) node {$y^1$};
\draw[anchor = east] (x21) node {$x_2^1$};
\draw[anchor = east] (z1) node {$z^1$};
\draw[anchor = north] (w11) node {$w_1^1$};
\draw[anchor = north] (w3n1) node {$w_{3n}^1$};

\draw (2.75,0.7) circle [x radius=3.9, y radius=1.6, rotate=90];
\node at (2.75,4) (n) {$H_n^1$};

\path (5,0) coordinate (x12);
\path (6.5,0) coordinate (y2);
\path (6.5,1.5) coordinate (x22);
\path (6.5,3) coordinate (z2);
\path (6,-1.5) coordinate (w12);
\path (7,-1.5) coordinate (w3n2);

\draw (y1) -- (x12);
\draw (x12) -- (y2) -- (x22) -- (z2);
\draw (y2) -- (w12);
\draw (y2) -- (w3n2);

\draw (x12) [fill=black] circle (\vr);
\draw (y2) [fill=black] circle (\vr);
\draw (x22) [fill=black] circle (\vr);
\draw (z2) [fill=white] circle (\vr);
\draw (w12) [fill=white] circle (\vr);
\node at (6.5,-1.5) (n) {$\dots$};
\draw (w3n2) [fill=white] circle (\vr);

\draw[anchor = north] (x12) node {$x_1^2$};
\draw[anchor = south east] (y2) node {$y^2$};
\draw[anchor = east] (x22) node {$x_2^2$};
\draw[anchor = east] (z2) node {$z^2$};
\draw[anchor = north] (w12) node {$w_1^2$};
\draw[anchor = north] (w3n2) node {$w_{3n}^2$};

\draw (6.25,0.7) circle [x radius=3.9, y radius=1.6, rotate=90];
\node at (6.25,4) (n) {$H_n^2$};

\path (8.5,0) coordinate (x13);
\draw (y2) -- (x13);
\draw (x13) [fill=black] circle (\vr);
\node at (9.5,0) (n) {$\dots$};
\draw[anchor = north] (x13) node {$x_1^3$};

\path (10.5,0) coordinate (x1n);
\path (12,0) coordinate (yn);
\path (12,1.5) coordinate (x2n);
\path (12,3) coordinate (zn);
\path (11.5,-1.5) coordinate (w1n);
\path (12.5,-1.5) coordinate (w3nn);

\draw (yn) -- (x1n);
\draw (x1n) -- (yn) -- (x2n) -- (zn);
\draw (yn) -- (w1n);
\draw (yn) -- (w3nn);

\draw (x1n) [fill=black] circle (\vr);
\draw (yn) [fill=black] circle (\vr);
\draw (x2n) [fill=black] circle (\vr);
\draw (zn) [fill=white] circle (\vr);
\draw (w1n) [fill=white] circle (\vr);
\node at (12,-1.5) (n) {$\dots$};
\draw (w3nn) [fill=white] circle (\vr);

\draw[anchor = north] (x1n) node {$x_1^n$};
\draw[anchor = south east] (yn) node {$y^n$};
\draw[anchor = east] (x2n) node {$x_2^n$};
\draw[anchor = east] (zn) node {$z^n$};
\draw[anchor = north] (w1n) node {$w_1^n$};
\draw[anchor = north] (w3nn) node {$w_{3n}^n$};

\draw (11.75,0.7) circle [x radius=3.9, y radius=1.6, rotate=90];
\node at (11.75,4) (n) {$H_n^n$};

\path (14,0) coordinate (x);
\path (15.5,0) coordinate (v);
\draw (yn) -- (x) -- (v);
\draw (x) [fill=black] circle (\vr);
\draw (v) [fill=white] circle (\vr);
\draw[anchor = north] (x) node {$x_1^{n+1}$};
\draw[anchor = north] (v) node {$v$};

\end{tikzpicture}
\end{center}
\caption{The infinite family of trees $T_n$ from Theorem \ref{S-m-Infinite}.}
\label{S-m-Infinite-Slika}
\end{figure}

\noindent {\textbf{Proof.}}
For an integer $n\geq 1$, let $H_n^i$, $i\in[n]$, be isomorphic graphs with 
$$V(H_n^i)=\{x_1^i,x_2^i,y^i,z^i\}\cup\{w_1^i,\dots,w_{3n}^i\}$$
and 
$$E(H_n^i)=\{x_1^iy^i, y^ix_2^i, x_2^iz^i\}\cup\{y^iw_j^i:j\in[3n]\}.$$
We define graph $T_n$ shown in Figure \ref{S-m-Infinite-Slika} by using graphs $H_n^i$, $i\in[n]$, as follows:
$$V(T_n)=\{u,v,x_1^{n+1}\}\cup \bigcup_{i=1}^{n}V(H_n^i)$$
and
$$E(T_n)=\{ux_1^1,x_1^{n+1}v,y^ix_1^{i+1}:i\in[n]\}\cup \bigcup_{i=1}^{n}E(H_n^i).$$

Notice that $d_{T_n}^s(x_1^1)=d_{T_n}^s(x_1^{n+1})=3n+3$, $d_{T_n}^s(x_1^i)=6n+5$, $i\in [2,n]$, $d_{T_n}^s(x_2^i)=3n+3$, $i\in [n]$, $d_{T_n}^s(y^1)=d_{T_n}^s(y^n)=3n+4$, $d_{T_n}^s(y^i)=3n+5$, $i\in [2,n-1]$, $d_{T_n}^s(w_j^i)=4$, $j\in [3n],\ i\in [n]$, and that the star degree of all other vertices of $T_n$ is at most $2$. Meaning that $m_s(T_n)=3n+1$. 


Next, we prove that $S_b(T_n)=2n+2$. To prove that $S_b(T_n)\geq 2n+2$ we construct a star b-coloring $c:V(T_n)\rightarrow[2n+2]$ in the following way:  
\begin{align*}
c(u)&=c(v)=c(y^i)=c(z^i)=2n+2, &&i\in [n],\\
c(x_1^i)&=2i-1, &&i\in [n+1],\\
c(x_2^i)&=2i, &&i\in [n].\\
\end{align*}
Additionally, for every $i\in [n]$ color vertices $w_j^i$, $j\in [2n-2]$, with pairwise distinct colors from the color set $[2n+1]-\{c(x_1^i),c(x_2^i),c(x_1^{i+1})\}$, and every vertex $w_j^i$, $j\in [2n-1,3n]$ with any color from $[2n+1]-\{c(x_1^i),c(x_2^i),c(x_1^{i+1})\}$.

It is easy to observe that every vertex $x_1^i$, $i\in[n+1]$, as well as every vertex $x_2^i$, $i\in [n]$, is a star b-vertex. Moreover, $y^1$ is also a star b-vertex. Thus, $S_b(T_n)\geq 2n+2$.

To show that $S_b(T_n)\leq 2n+2$, we need to check that there does not exist a star b-coloring of $T_n$ using at least $2n+3$ colors. Notice that if there exists a star coloring such that vertices $y^i$ and $y^{i+1}$, $i\in [n-1]$, receive distinct colors, then vertex $x_1^{i+1}$ cannot be a star b-vertex of $c$. So, for every possible star b-vertex $y^i$, $i>1$, of color different than $y^{i-1}$, we have one b-vertex less, namely $x_1^{i-1}$. This means that it is not possible to obtain a star b-coloring of $T_n$ using $2n+3$ colors. Now, $m_s(T_n)-S_b(T_n)=n-1$ tends to infinity as tends to infinity as $n$ increases.
\qed\bigskip

The same family of trees as in above proof also yields an arbitrary difference between $S_b(G)$ and $S(G)$.

\begin{theorem}\label{infinite}
There exists an infinite family of trees $T_1, T_2, \dots$ such that $(S_b(T_n)-S(T_n))\rightarrow\infty$ as $n\rightarrow\infty$. 
\end{theorem}

\noindent {\textbf{Proof.}}
Let $T_1, T_2, \dots$ be a family of trees from the proof of Theorem \ref{S-m-Infinite}. We have already shown that $S_b(T_n)=2n+2$. In order to prove that $S(T_n)\leq 3$ let $k=\left\lfloor \frac{n}{2}\right\rfloor$ and we construct a star coloring $c:V(T_n)\to [3]$ as follows:
\begin{align*}
c(u)&=3,\\
c(x_1^i)&=1, &&i\in [n+1],\\
c(x_2^i)&=c(w_j^i)=1, &&i\in [n],j\in[3n],\\
c(y^{2i-1})&=c(z^{2i})=2, &&i\in [k],\\
c(y^{2i})&=c(z^{2i-1})=3, &&i\in [k],\\
\end{align*}
and $c(v)=2$ for even $n$ or $c(y^n)=2$ and $c(v)=3$ for odd $n$. Clearly, the obtained coloring $c$ is a proper coloring of $T_n$ that does not induce a bi-chromatic path on 4 vertices. Therefore, $S(T_n)\leq 3$ and $(S_b(T_n)-S(T_n))\rightarrow\infty$ as $n\rightarrow\infty$. 
\qed\bigskip

Next, we obtain some results on star b-coloring of joins. 
Let $\sqcup $ denote the disjoint union. Recall that the join of graphs $G$ and $H$ is the graph $G\vee H$ with the vertex set $V(G\vee H)=V(G)\sqcup V(H)$ and the edge set $E(G\vee H)=E(G)\sqcup E(H)\sqcup \{uv:u\in V(G)\wedge v\in V(H)\}$.
The following result can be obtained analogously as the acyclic b-chromatic number of joins in \cite{AnCP}. The proof is the same only that we need to replace acyclic b-colorings with star b-colorings and acyclic b-chromatic numbers with star b-chromatic numbers whenever they appear in the proof. Therefore we omit this proof.
 
\begin{theorem}\label{join}
For two non-complete graphs $G$ and $H$ we have 
$$S_b(G\vee H)=\max\{S_b(G)+n_H,S_b(H)+n_G\}.$$
If $H\cong K_q$, then $S_b(G\vee H)=S_b(G)+q$.
\end{theorem}

Note that complete bipartite graph $K_{m,n}=\overline{K}_n\vee\overline{K}_m$, wheel $W_n=K_1\vee C_{n-1}$, fan $F_n=K_1\vee P_{n-1}$ and complete split graph $K_n\vee\overline{K}_m$ can all be obtained as joins of two graphs. Thus the next corollary is a direct consequence of Theorem \ref{join}. In particular, second and third item can be further applied by Theorems \ref{poti} and \ref{cikli}, respectively, and the last item follows from Proposition \ref{order}.

\begin{corollary}\label{wheel}
For every positive integers $k,\ell,m,n$, where $k,\ell\geq 4$,  we have 
\begin{itemize}
\item $S_b(K_{n,m})=1+\max\{n,m\}$; 
\item $S_b(W_k)=S_b(C_{k-1})+1$; 
\item $S_b(F_{\ell})=S_b(P_{\ell-1})+1$; 
\item $S_b(K_n\vee\overline{K}_m)=n+1$. 
\end{itemize}
\end{corollary}

Using Corollary \ref{wheel} we infer that the difference between $S_b(G)$ and $\varphi(G)$ can be arbitrarily large. (The familly of trees from Theorem \ref{DeltaSquaredExtremal} is also an example for the next result.)

\begin{corollary}\label{wheelInfinite}
There exists an infinite family of graphs $G_1, G_2, \dots$ such that $(S_b(G_n)-\varphi(G_n))\rightarrow\infty$ as $n\rightarrow\infty$. 
\end{corollary}

\noindent {\textbf{Proof.}}
Let $G_n=K_{n,n}$ for any $n\geq 1$. It is easy to observe that $\varphi(G_n)=2$, while $S_b(G_n)=n+1$ follows from Theorem \ref{wheel}. Thus, the result follows.
\qed\bigskip



\section{Open questions}

We conclude this paper with several open questions. The first concerns the star degree of a vertex, which plays a crucial role in determining $m_s(G)$, an upper bound for $S_b(G)$. Therefore, an efficient method for computing $d_G^s(v)$ would be valuable. While we solved this for certain special cases in Theorems \ref{girth} and \ref{matching}, the general case remains open.

\begin{question}
	Can $d_G^s(v)$ be computed in polynomial time for a general graph $G$, or is this problem NP-hard?
\end{question}

The computational complexity of the star b-chromatic number remains an open question. We attempted to approach it from several perspectives but were not successful.

\begin{problem}
	What is the computational complexity of determining $S_b(G)$?
\end{problem}

We assume that the problem itself does not have a polynomial-time solution. Therefore, it would be valuable to identify some graph classes for which $S_b(G)$ can be computed in polynomial time.

\begin{problem}
	Find graph classes where $S_b(G)$ can be computed in polynomial time.
\end{problem}

By Theorem \ref{girth} we can compute $d_T^s(v)$ for any tree $T$. However, Theorems \ref{DeltaSquaredExtremal}, \ref{poti} and \ref{S-m-Infinite} imply that $S_b(T)$ is much more complex than $\varphi(T)$ and $A_b(T)$ (see \cite{IrMa} and \cite{AnCP}, respectively). This leads to the following question. 

\begin{question}
	Determine $S_b(T)$ for any tree $T$.
\end{question}


\end{document}